\documentclass[runningheads]{cl2emult}

\usepackage{makeidx}  
\usepackage{subeqnar} 
\usepackage{multicol} 
\usepackage{cropmark} 
\usepackage{math}     
\makeindex            

\usepackage{amsmath,amsfonts}   
\catcode`\@=11

\thinlines
\newskip\Einheit \Einheit=0.6cm
\newcount\xcoord \newcount\ycoord
\newdimen\xdim \newdimen\ydim \newdimen\PfadD@cke \newdimen\Pfadd@cke
\def\PfadDicke#1{\PfadD@cke#1 \divide\PfadD@cke by2 \Pfadd@cke\PfadD@cke \multiply\PfadD@cke by2}
\long\def\LOOP#1\REPEAT{\def\BODY{#1}\ITERATE}
\def\ITERATE{\BODY \let\next\ITERATE \else\let\next\relax\fi \next}
\let\REPEAT=\fi
\def\Punkt{\hbox{\raise-2pt\hbox to0pt{\hss\scriptsize$\bullet$\hss}}}
\def\DuennPunkt(#1,#2){\unskip
  \raise#2 \Einheit\hbox to0pt{\hskip#1 \Einheit
          \raise-2.5pt\hbox to0pt{\hss\normalsize$\bullet$\hss}\hss}}
\def\NormalPunkt(#1,#2){\unskip
  \raise#2 \Einheit\hbox to0pt{\hskip#1 \Einheit
          \raise-3pt\hbox to0pt{\hss\large$\bullet$\hss}\hss}}
\def\DickPunkt(#1,#2){\unskip
  \raise#2 \Einheit\hbox to0pt{\hskip#1 \Einheit
          \raise-4pt\hbox to0pt{\hss\Large$\bullet$\hss}\hss}}
\def\Kreis(#1,#2){\unskip
  \raise#2 \Einheit\hbox to0pt{\hskip#1 \Einheit
          \raise-4pt\hbox to0pt{\hss\Large$\circ$\hss}\hss}}
\def\Diagonale(#1,#2)#3{\unskip\leavevmode
  \xcoord#1\relax \ycoord#2\relax
      \raise\ycoord \Einheit\hbox to0pt{\hskip\xcoord \Einheit
         \unitlength\Einheit
         \line(1,1){#3}\hss}}
\def\AntiDiagonale(#1,#2)#3{\unskip\leavevmode
  \xcoord#1\relax \ycoord#2\relax \advance\xcoord by -0.05\relax
      \raise\ycoord \Einheit\hbox to0pt{\hskip\xcoord \Einheit
         \unitlength\Einheit
         \line(1,-1){#3}\hss}}
\def\Pfad(#1,#2),#3\endPfad{\unskip\leavevmode
  \xcoord#1 \ycoord#2 \thicklines\ZeichnePfad#3\endPfad\thinlines}
\def\ZeichnePfad#1{\ifx#1\endPfad\let\next\relax
  \else\let\next\ZeichnePfad
    \ifnum#1=1
      \raise\ycoord \Einheit\hbox to0pt{\hskip\xcoord \Einheit
         \vrule height\Pfadd@cke width1 \Einheit depth\Pfadd@cke\hss}%
      \advance\xcoord by 1
    \else\ifnum#1=2
      \raise\ycoord \Einheit\hbox to0pt{\hskip\xcoord \Einheit
        \hbox{\hskip-\PfadD@cke\vrule height1 \Einheit width\PfadD@cke depth0pt}\hss}%
      \advance\ycoord by 1
    \else\ifnum#1=3
      \raise\ycoord \Einheit\hbox to0pt{\hskip\xcoord \Einheit
         \unitlength\Einheit
         \line(1,1){1}\hss}
      \advance\xcoord by 1
      \advance\ycoord by 1
    \else\ifnum#1=4
      \raise\ycoord \Einheit\hbox to0pt{\hskip\xcoord \Einheit
         \unitlength\Einheit
         \line(1,-1){1}\hss}
      \advance\xcoord by 1
      \advance\ycoord by -1
    \else\ifnum#1=5
      \advance\xcoord by -1
      \raise\ycoord \Einheit\hbox to0pt{\hskip\xcoord \Einheit
         \vrule height\Pfadd@cke width1 \Einheit depth\Pfadd@cke\hss}%
    \else\ifnum#1=6
      \advance\ycoord by -1
      \raise\ycoord \Einheit\hbox to0pt{\hskip\xcoord \Einheit
        \hbox{\hskip-\PfadD@cke\vrule height1 \Einheit width\PfadD@cke depth0pt}\hss}%
    \else\ifnum#1=7
      \advance\xcoord by -1
      \advance\ycoord by -1
      \raise\ycoord \Einheit\hbox to0pt{\hskip\xcoord \Einheit
         \Line@(1,1){1}\hss}
    \else\ifnum#1=8
      \advance\xcoord by -1
      \advance\ycoord by +1
      \raise\ycoord \Einheit\hbox to0pt{\hskip\xcoord \Einheit
         \Line@(1,-1){1}\hss}
    \fi\fi\fi\fi
    \fi\fi\fi\fi
  \fi\next}
\def\hSSchritt{\leavevmode\raise-.4pt\hbox to0pt{\hss.\hss}\hskip.2\Einheit
  \raise-.4pt\hbox to0pt{\hss.\hss}\hskip.2\Einheit
  \raise-.4pt\hbox to0pt{\hss.\hss}\hskip.2\Einheit
  \raise-.4pt\hbox to0pt{\hss.\hss}\hskip.2\Einheit
  \raise-.4pt\hbox to0pt{\hss.\hss}\hskip.2\Einheit}
\def\vSSchritt{\vbox{\baselineskip.2\Einheit\lineskiplimit0pt
\hbox{.}\hbox{.}\hbox{.}\hbox{.}\hbox{.}}}
\def\DSSchritt{\leavevmode\raise-.4pt\hbox to0pt{%
  \hbox to0pt{\hss.\hss}\hskip.2\Einheit
  \raise.2\Einheit\hbox to0pt{\hss.\hss}\hskip.2\Einheit
  \raise.4\Einheit\hbox to0pt{\hss.\hss}\hskip.2\Einheit
  \raise.6\Einheit\hbox to0pt{\hss.\hss}\hskip.2\Einheit
  \raise.8\Einheit\hbox to0pt{\hss.\hss}\hss}}
\def\dSSchritt{\leavevmode\raise-.4pt\hbox to0pt{%
  \hbox to0pt{\hss.\hss}\hskip.2\Einheit
  \raise-.2\Einheit\hbox to0pt{\hss.\hss}\hskip.2\Einheit
  \raise-.4\Einheit\hbox to0pt{\hss.\hss}\hskip.2\Einheit
  \raise-.6\Einheit\hbox to0pt{\hss.\hss}\hskip.2\Einheit
  \raise-.8\Einheit\hbox to0pt{\hss.\hss}\hss}}
\def\SPfad(#1,#2),#3\endSPfad{\unskip\leavevmode
  \xcoord#1 \ycoord#2 \ZeichneSPfad#3\endSPfad}
\def\ZeichneSPfad#1{\ifx#1\endSPfad\let\next\relax
  \else\let\next\ZeichneSPfad
    \ifnum#1=1
      \raise\ycoord \Einheit\hbox to0pt{\hskip\xcoord \Einheit
         \hSSchritt\hss}%
      \advance\xcoord by 1
    \else\ifnum#1=2
      \raise\ycoord \Einheit\hbox to0pt{\hskip\xcoord \Einheit
        \hbox{\hskip-2pt \vSSchritt}\hss}%
      \advance\ycoord by 1
    \else\ifnum#1=3
      \raise\ycoord \Einheit\hbox to0pt{\hskip\xcoord \Einheit
         \DSSchritt\hss}
      \advance\xcoord by 1
      \advance\ycoord by 1
    \else\ifnum#1=4
      \raise\ycoord \Einheit\hbox to0pt{\hskip\xcoord \Einheit
         \dSSchritt\hss}
      \advance\xcoord by 1
      \advance\ycoord by -1
    \else\ifnum#1=5
      \advance\xcoord by -1
      \raise\ycoord \Einheit\hbox to0pt{\hskip\xcoord \Einheit
         \hSSchritt\hss}%
    \else\ifnum#1=6
      \advance\ycoord by -1
      \raise\ycoord \Einheit\hbox to0pt{\hskip\xcoord \Einheit
        \hbox{\hskip-2pt \vSSchritt}\hss}%
    \else\ifnum#1=7
      \advance\xcoord by -1
      \advance\ycoord by -1
      \raise\ycoord \Einheit\hbox to0pt{\hskip\xcoord \Einheit
         \DSSchritt\hss}
    \else\ifnum#1=8
      \advance\xcoord by -1
      \advance\ycoord by 1
      \raise\ycoord \Einheit\hbox to0pt{\hskip\xcoord \Einheit
         \dSSchritt\hss}
    \fi\fi\fi\fi
    \fi\fi\fi\fi
  \fi\next}
\def\Koordinatenachsen(#1,#2){\unskip
 \hbox to0pt{\hskip-.5pt\vrule height#2 \Einheit width.5pt depth1 \Einheit}%
 \hbox to0pt{\hskip-1 \Einheit \xcoord#1 \advance\xcoord by1
    \vrule height0.25pt width\xcoord \Einheit depth0.25pt\hss}}
\def\Koordinatenachsen(#1,#2)(#3,#4){\unskip
 \hbox to0pt{\hskip-.5pt \ycoord-#4 \advance\ycoord by1
    \vrule height#2 \Einheit width.5pt depth\ycoord \Einheit}%
 \hbox to0pt{\hskip-1 \Einheit \hskip#3\Einheit 
    \xcoord#1 \advance\xcoord by1 \advance\xcoord by-#3 
    \vrule height0.25pt width\xcoord \Einheit depth0.25pt\hss}}
\def\Gitter(#1,#2){\unskip \xcoord0 \ycoord0 \leavevmode
  \LOOP\ifnum\ycoord<#2
    \loop\ifnum\xcoord<#1
      \raise\ycoord \Einheit\hbox to0pt{\hskip\xcoord \Einheit\Punkt\hss}%
      \advance\xcoord by1
    \repeat
    \xcoord0
    \advance\ycoord by1
  \REPEAT}
\def\Gitter(#1,#2)(#3,#4){\unskip \xcoord#3 \ycoord#4 \leavevmode
  \LOOP\ifnum\ycoord<#2
    \loop\ifnum\xcoord<#1
      \raise\ycoord \Einheit\hbox to0pt{\hskip\xcoord \Einheit\Punkt\hss}%
      \advance\xcoord by1
    \repeat
    \xcoord#3
    \advance\ycoord by1
  \REPEAT}
\def\Label#1#2(#3,#4){\unskip \xdim#3 \Einheit \ydim#4 \Einheit
  \def\lo{\advance\xdim by-.5 \Einheit \advance\ydim by.5 \Einheit}%
  \def\llo{\advance\xdim by-.25cm \advance\ydim by.5 \Einheit}%
  \def\loo{\advance\xdim by-.5 \Einheit \advance\ydim by.25cm}%
  \def\o{\advance\ydim by.25cm}%
  \def\ro{\advance\xdim by.5 \Einheit \advance\ydim by.5 \Einheit}%
  \def\rro{\advance\xdim by.25cm \advance\ydim by.5 \Einheit}%
  \def\roo{\advance\xdim by.5 \Einheit \advance\ydim by.25cm}%
  \def\l{\advance\xdim by-.30cm}%
  \def\r{\advance\xdim by.30cm}%
  \def\lu{\advance\xdim by-.5 \Einheit \advance\ydim by-.6 \Einheit}%
  \def\llu{\advance\xdim by-.25cm \advance\ydim by-.6 \Einheit}%
  \def\luu{\advance\xdim by-.5 \Einheit \advance\ydim by-.30cm}%
  \def\u{\advance\ydim by-.30cm}%
  \def\ru{\advance\xdim by.5 \Einheit \advance\ydim by-.6 \Einheit}%
  \def\rru{\advance\xdim by.25cm \advance\ydim by-.6 \Einheit}%
  \def\ruu{\advance\xdim by.5 \Einheit \advance\ydim by-.30cm}%
  #1\raise\ydim\hbox to0pt{\hskip\xdim
     \vbox to0pt{\vss\hbox to0pt{\hss$#2$\hss}\vss}\hss}%
}

\newcounter{saveeqn}
\newcommand{\alphaeqn}{\setcounter{saveeqn}{\value{equation}}%
\setcounter{equation}{0}%
\if@numart
\global\def\theequation{\arabic{saveeqn}\alph{equation}}%
\else
\global\def\theequation{\thechapter.\arabic{saveeqn}\alph{equation}}%
\fi}
\newcommand{\reseteqn}{\setcounter{equation}{\value{saveeqn}}%
\if@numart
\global\def\theequation{\arabic{equation}}%
\else
\global\def\theequation{\thechapter.\arabic{equation}}%
\fi}

\catcode`\@=12

\def\al{\alpha}
\def\ep{\varepsilon}
\def\si{\sigma}
\def\ph{\varphi}
\def\P{{\mathcal P}}
\def\S{{\mathfrak S}}

\def\TA{\text{\it TA}}
\def\GF{\operatorname{GF}}
\def\NE{\operatorname{NE}}
\def\sgn{\operatorname{sgn}}
\def\lK{\left(}
\def\rK{\right)}

\begin{document}
\title*{A determinantal formula for the Hilbert series of one-sided ladder
determinantal rings}
\toctitle{A determinantal formula for the Hilbert series of one-sided ladder
determinantal rings}
%
%
\titlerunning{The Hilbert series of ladder determinantal rings}
%
\author{C.~Krattenthaler$^\ast$ and M.~Rubey\thanks{Research partially
supported by the Austrian
Science Foundation FWF, grant P12094-MAT.} 
}
\authorrunning{C. Krattenthaler and M. Rubey}
%
%
\institute{
Institut f\"ur Mathematik der Universit\"at Wien,\\
Strudlhofgasse 4, A-1090 Wien, Austria.\\
e-mail: KRATT@Ap.Univie.Ac.At, a9104910@unet.univie.ac.at\\
WWW: \tt http://www.mat.univie.ac.at/\~{}kratt
}


\maketitle              

\bigskip
\centerline{\small\it Dedicated to Shreeram Abhyankar}

\begin{abstract}
We give a formula that expresses the Hilbert series of
one-sided ladder determinantal rings, up to a trivial factor, in form
of a determinant.  This allows the convenient computation of these
Hilbert series.  The formula follows from a determinantal formula for
a generating function for families of nonintersecting lattice paths
that stay inside a one-sided ladder-shaped region, in which the paths
are counted with respect to turns.
\end{abstract}

\section{Introduction}
\label{Sec:1}

Work of Abhyankar and Kulkarni
\cite{AbhyAB,AbKuAC,KulkAD,KulkAC}, Bruns, Conca, Herzog, and
Trung \cite{BrHeAA,ConcAB,CoHeAA,HeTrAA} showed that the
computation of the Hilbert series of ladder determinantal rings (see
Section~\ref{Sec:2} for precise definitions and background) boils down to
counting families of $n$ nonintersecting lattice paths with a given
total number of turns in a certain ladder-shaped region. Thus, this
raises the question of establishing an explicit formula for the number
of these families of nonintersecting lattice paths.

In the case that there is no ladder restriction, Abhyankar
\cite[(20.14.4)]{AbhyAB} has found a determinantal formula for the Hilbert
series (actually not just one, but a great number of them). 
As was made explicit in \cite{CoHeAA,GhorAC,KulkAC,ModaAA}, 
he thereby solved the aforementioned counting problem in the
case of {\it no} ladder restriction. For direct proofs of the corresponding
counting formula see \cite{KratBE,ModaAA}.  In the case of
{\it one-sided\/} ladders, Kulkarni \cite{KulkAD} established an explicit
solution to the counting problem for $n=1$ (i.e., if there is just one
path; this corresponds to considering one-sided ladder determinantal
rings defined by $2\times 2$ minors). For arbitrary $n$, a
determinantal formula for the number of families of $n$
nonintersecting lattice paths in a one-sided ladder, where the
starting and end points of the paths are successive, was given by the
first author and Prohaska \cite{KrPrAA} (this corresponds to one-sided
ladder determinantal rings defined by $(n+1)\times (n+1)$ minors),
thereby proving a conjecture by Conca and Herzog \cite[last
paragraph]{CoHeAA}. Finally, Ghorpade \cite{GhorAF} has recently proposed a
solution to the counting problem with more general starting and end
points of the paths, even in the case of {\it two-sided\/} ladders (this
corresponds to {\it two-sided\/} ladder determinantal rings {\it
cogenerated\/} by a
given minor). This solution is based on an explicit formula for the
counting problem for one path (i.e., $n=1$), which is then summed over
a large set of indices with complicated dependencies. Thus, this
solution cannot be regarded as equally satisfying as the determinantal
formula of Abhyankar and the determinantal formula of the first author and
Prohaska, which are, however, only formulas in the case of a trivial
ladder and in the case of a one-sided ladder, respectively.

The purpose of this paper is to provide a determinantal formula for
the case of {\it one-sided\/} ladders where the starting and end points are
more general than in \cite{KrPrAA} (see Corollary~\ref{cor:3}; this corresponds
to {\it one-sided\/} ladder determinantal rings {\it cogenerated\/} by a given
minor). This formula must be considered as superior to the
aforementioned one by Ghorpade \cite{GhorAF} in this case (i.e., the
case of one- instead of two-sided ladders). It specializes directly to
Abhyankar's formula \cite[(20.14.4), $L=2$, $k=2$, with
$F^{(22)}(m,p,a,V)$ defined on p.~50]{AbhyAB} in the case of no ladder
restriction. On the other hand, if starting and end points are
successive, then it does not specialize to the formula in
\cite{KrPrAA}. (As already mentioned in Section~7 of \cite{KrPrAA},
it seems that the formula in \cite{KrPrAA} cannot be extended in any
direction.)

The entries in the determinant in our formula \eqref{e3.3}, respectively
\eqref{e3.4}, are given by certain generating functions for two-rowed arrays,
which are easy to compute as we show in Section~\ref{Sec:5}. (The concept of
two-rowed arrays was introduced in \cite{KratAF,KrMoAC} and
developed to full power in \cite{KratAP,KratBE}. Also the
proof of the main theorem in \cite{KrPrAA} depended heavily on
two-rowed arrays.)

In the next section we recall the basic setup. In particular, we
define ladder determinantal rings and state, in Theorem~\ref{thm:1}, the
connection between the Hilbert series of such rings and the
enumeration of nonintersecting lattice paths with respect to turns.
Our main result, the determinantal formula for the Hilbert series of
one-sided ladder determinantal rings cogenerated by a given fixed
minor, is stated in Corollary~\ref{cor:3} in Section~\ref{Sec:3}. It
follows from a determinantal formula for counting nonintersecting
lattice paths in a one-sided ladder with respect to turns, where the
starting and end points are allowed to be even more general than is
needed for our main result. This counting formula is stated in
Theorem~\ref{thm:2}, and it is proved in Section~\ref{Sec:4}. 
In Section~\ref{Sec:5} we show how to compute the
generating functions for two-rowed arrays that appear in the
determinant of our formula.

\section{Ladder determinantal rings and the enumeration of
nonintersecting lattice paths with respect to turns}
\label{Sec:2}

Let $X=(X_{i,j})_{0\le i\le b,\ 0\le j\le a}$ be a $(b+1)\times (a+1)$
matrix of indeterminates. Let $=(Y_{i,j})_{0\le i\le b,\ 0\le j\le a}$
be another $(b+1)\times (a+1)$ matrix with the property that
$Y_{i,j}=X_{i,j}$ or 0, and if $Y_{i,j}=X_{i,j}$ and
$Y_{i'j'}=X_{i'j'}$, where $i\le i'$ and $j\le j'$, then
$Y_{s,t}=X_{s,t}$ for all $s,t$ with $i\le s\le i'$ and $j\le t\le
j'$. An example for such a matrix $Y$, with $b=15$ and $a=13$ is
displayed in Figure~1.  (Note that there could be 0's in the
bottom-right corner of the matrix also.)  Such a ``submatrix" $Y$ of
$X$ is called a {\it ladder}. This terminology is motivated by the
identification of such a matrix $Y$ with the set of all points
$(j,b-i)$ in the plane for which $Y_{i,j}=X_{i,j}$. For example, the
set of all such points for the special matrix in Figure~1 is shown in
Figure~2. (It should be apparent from comparison of
Figures~1 and 2
that the reason for taking $(j,b-i)$ instead of $(i,j)$ is to take
care of the difference in ``orientation" of row and column indexing of
a matrix versus coordinates in the plane.)  In general, this set of
points looks like a (two-sided) ladder-shaped region. If, on the other
hand, we have either $Y_{0,0}=X_{0,0}$ or $Y_{b,a}=X_{b,a}$ then we
call $Y$ a {\it one-sided\/} ladder. In the first case we call $Y$ a
{\it lower ladder}, in the second an {\it upper ladder}.  Thus, the
matrix in Figure~1 is an upper ladder.

\begin{figure}[h] \label{fig1}
$$\Einheit.6cm
\left(
\hbox{\hskip.4cm}
\Label\o{\hbox{\scriptsize$X\!\!_{1\!5,0}$}}(0,-8)
\Label\o{\hbox{\scriptsize$X\!\!_{1\!5,1}$}}(1,-8)
\Label\o{\hbox{\scriptsize$X\!\!_{1\!5,2}$}}(2,-8)
\Label\o{\hbox{\scriptsize$X\!\!_{1\!5,3}$}}(3,-8)
\Label\o{\hbox{\scriptsize$X\!\!_{1\!5,4}$}}(4,-8)
\Label\o{\hbox{\scriptsize$X\!\!_{1\!5,5}$}}(5,-8)
\Label\o{\hbox{\scriptsize$X\!\!_{1\!5,6}$}}(6,-8)
\Label\o{\hbox{\scriptsize$X\!\!_{1\!5,7}$}}(7,-8)
\Label\o{\hbox{\scriptsize$X\!\!_{1\!5,8}$}}(8,-8)
\Label\o{\hbox{\scriptsize$X\!\!_{1\!5,9}$}}(9,-8)
\Label\o{\hbox{\scriptsize$X\!\!_{1\!5,1\!0}$}}(10,-8)
\Label\o{\hbox{\scriptsize$X\!\!_{1\!5,1\!1}$}}(11,-8)
\Label\o{\hbox{\scriptsize$X\!\!_{1\!5,1\!2}$}}(12,-8)
\Label\o{\hbox{\scriptsize$X\!\!_{1\!5,1\!3}$}}(13,-8)
\Label\o{\hbox{\scriptsize$X\!\!_{1\!4,0}$}}(0,-7)
\Label\o{\hbox{\scriptsize$X\!\!_{1\!4,1}$}}(1,-7)
\Label\o{\hbox{\scriptsize$X\!\!_{1\!4,2}$}}(2,-7)
\Label\o{\hbox{\scriptsize$X\!\!_{1\!4,3}$}}(3,-7)
\Label\o{\hbox{\scriptsize$X\!\!_{1\!4,4}$}}(4,-7)
\Label\o{\hbox{\scriptsize$X\!\!_{1\!4,5}$}}(5,-7)
\Label\o{\hbox{\scriptsize$X\!\!_{1\!4,6}$}}(6,-7)
\Label\o{\hbox{\scriptsize$X\!\!_{1\!4,7}$}}(7,-7)
\Label\o{\hbox{\scriptsize$X\!\!_{1\!4,8}$}}(8,-7)
\Label\o{\hbox{\scriptsize$X\!\!_{1\!4,9}$}}(9,-7)
\Label\o{\hbox{\scriptsize$X\!\!_{1\!4,1\!0}$}}(10,-7)
\Label\o{\hbox{\scriptsize$X\!\!_{1\!4,1\!1}$}}(11,-7)
\Label\o{\hbox{\scriptsize$X\!\!_{1\!4,1\!2}$}}(12,-7)
\Label\o{\hbox{\scriptsize$X\!\!_{1\!4,1\!3}$}}(13,-7)
\Label\o{\hbox{\scriptsize$X\!\!_{1\!3,0}$}}(0,-6)
\Label\o{\hbox{\scriptsize$X\!\!_{1\!3,1}$}}(1,-6)
\Label\o{\hbox{\scriptsize$X\!\!_{1\!3,2}$}}(2,-6)
\Label\o{\hbox{\scriptsize$X\!\!_{1\!3,3}$}}(3,-6)
\Label\o{\hbox{\scriptsize$X\!\!_{1\!3,4}$}}(4,-6)
\Label\o{\hbox{\scriptsize$X\!\!_{1\!3,5}$}}(5,-6)
\Label\o{\hbox{\scriptsize$X\!\!_{1\!3,6}$}}(6,-6)
\Label\o{\hbox{\scriptsize$X\!\!_{1\!3,7}$}}(7,-6)
\Label\o{\hbox{\scriptsize$X\!\!_{1\!3,8}$}}(8,-6)
\Label\o{\hbox{\scriptsize$X\!\!_{1\!3,9}$}}(9,-6)
\Label\o{\hbox{\scriptsize$X\!\!_{1\!3,1\!0}$}}(10,-6)
\Label\o{\hbox{\scriptsize$X\!\!_{1\!3,1\!1}$}}(11,-6)
\Label\o{\hbox{\scriptsize$X\!\!_{1\!3,1\!2}$}}(12,-6)
\Label\o{\hbox{\scriptsize$X\!\!_{1\!3,1\!3}$}}(13,-6)
\Label\o{\hbox{\scriptsize$X\!\!_{1\!2,0}$}}(0,-5)
\Label\o{\hbox{\scriptsize$X\!\!_{1\!2,1}$}}(1,-5)
\Label\o{\hbox{\scriptsize$X\!\!_{1\!2,2}$}}(2,-5)
\Label\o{\hbox{\scriptsize$X\!\!_{1\!2,3}$}}(3,-5)
\Label\o{\hbox{\scriptsize$X\!\!_{1\!2,4}$}}(4,-5)
\Label\o{\hbox{\scriptsize$X\!\!_{1\!2,5}$}}(5,-5)
\Label\o{\hbox{\scriptsize$X\!\!_{1\!2,6}$}}(6,-5)
\Label\o{\hbox{\scriptsize$X\!\!_{1\!2,7}$}}(7,-5)
\Label\o{\hbox{\scriptsize$X\!\!_{1\!2,8}$}}(8,-5)
\Label\o{\hbox{\scriptsize$X\!\!_{1\!2,9}$}}(9,-5)
\Label\o{\hbox{\scriptsize$X\!\!_{1\!2,1\!0}$}}(10,-5)
\Label\o{\hbox{\scriptsize$X\!\!_{1\!2,1\!1}$}}(11,-5)
\Label\o{\hbox{\scriptsize$X\!\!_{1\!2,1\!2}$}}(12,-5)
\Label\o{\hbox{\scriptsize$X\!\!_{1\!2,1\!3}$}}(13,-5)
\Label\o{\hbox{\scriptsize$X\!\!_{1\!1,0}$}}(0,-4)
\Label\o{\hbox{\scriptsize$X\!\!_{1\!1,1}$}}(1,-4)
\Label\o{\hbox{\scriptsize$X\!\!_{1\!1,2}$}}(2,-4)
\Label\o{\hbox{\scriptsize$X\!\!_{1\!1,3}$}}(3,-4)
\Label\o{\hbox{\scriptsize$X\!\!_{1\!1,4}$}}(4,-4)
\Label\o{\hbox{\scriptsize$X\!\!_{1\!1,5}$}}(5,-4)
\Label\o{\hbox{\scriptsize$X\!\!_{1\!1,6}$}}(6,-4)
\Label\o{\hbox{\scriptsize$X\!\!_{1\!1,7}$}}(7,-4)
\Label\o{\hbox{\scriptsize$X\!\!_{1\!1,8}$}}(8,-4)
\Label\o{\hbox{\scriptsize$X\!\!_{1\!1,9}$}}(9,-4)
\Label\o{\hbox{\scriptsize$X\!\!_{1\!1,1\!0}$}}(10,-4)
\Label\o{\hbox{\scriptsize$X\!\!_{1\!1,1\!1}$}}(11,-4)
\Label\o{\hbox{\scriptsize$X\!\!_{1\!1,1\!2}$}}(12,-4)
\Label\o{\hbox{\scriptsize$X\!\!_{1\!1,1\!3}$}}(13,-4)
\Label\o{\hbox{\scriptsize$X\!\!_{1\!0,0}$}}(0,-3)
\Label\o{\hbox{\scriptsize$X\!\!_{1\!0,1}$}}(1,-3)
\Label\o{\hbox{\scriptsize$X\!\!_{1\!0,2}$}}(2,-3)
\Label\o{\hbox{\scriptsize$X\!\!_{1\!0,3}$}}(3,-3)
\Label\o{\hbox{\scriptsize$X\!\!_{1\!0,4}$}}(4,-3)
\Label\o{\hbox{\scriptsize$X\!\!_{1\!0,5}$}}(5,-3)
\Label\o{\hbox{\scriptsize$X\!\!_{1\!0,6}$}}(6,-3)
\Label\o{\hbox{\scriptsize$X\!\!_{1\!0,7}$}}(7,-3)
\Label\o{\hbox{\scriptsize$X\!\!_{1\!0,8}$}}(8,-3)
\Label\o{\hbox{\scriptsize$X\!\!_{1\!0,9}$}}(9,-3)
\Label\o{\hbox{\scriptsize$X\!\!_{1\!0,1\!0}$}}(10,-3)
\Label\o{\hbox{\scriptsize$X\!\!_{1\!0,1\!1}$}}(11,-3)
\Label\o{\hbox{\scriptsize$X\!\!_{1\!0,1\!2}$}}(12,-3)
\Label\o{\hbox{\scriptsize$X\!\!_{1\!0,1\!3}$}}(13,-3)
\Label\o{\hbox{\scriptsize$X\!\!_{9,0}$}}(0,-2)
\Label\o{\hbox{\scriptsize$X\!\!_{9,1}$}}(1,-2)
\Label\o{\hbox{\scriptsize$X\!\!_{9,2}$}}(2,-2)
\Label\o{\hbox{\scriptsize$X\!\!_{9,3}$}}(3,-2)
\Label\o{\hbox{\scriptsize$X\!\!_{9,4}$}}(4,-2)
\Label\o{\hbox{\scriptsize$X\!\!_{9,5}$}}(5,-2)
\Label\o{\hbox{\scriptsize$X\!\!_{9,6}$}}(6,-2)
\Label\o{\hbox{\scriptsize$X\!\!_{9,7}$}}(7,-2)
\Label\o{\hbox{\scriptsize$X\!\!_{9,8}$}}(8,-2)
\Label\o{\hbox{\scriptsize$X\!\!_{9,9}$}}(9,-2)
\Label\o{\hbox{\scriptsize$X\!\!_{9,1\!0}$}}(10,-2)
\Label\o{\hbox{\scriptsize$X\!\!_{9,1\!1}$}}(11,-2)
\Label\o{\hbox{\scriptsize$X\!\!_{9,1\!2}$}}(12,-2)
\Label\o{\hbox{\scriptsize$X\!\!_{9,1\!3}$}}(13,-2)
\Label\o{\hbox{\scriptsize$0$}}(0,-1)
\Label\o{\hbox{\scriptsize$0$}}(1,-1)
\Label\o{\hbox{\scriptsize$0$}}(2,-1)
\Label\o{\hbox{\scriptsize$0$}}(3,-1)
\Label\o{\hbox{\scriptsize$X\!\!_{8,4}$}}(4,-1)
\Label\o{\hbox{\scriptsize$X\!\!_{8,5}$}}(5,-1)
\Label\o{\hbox{\scriptsize$X\!\!_{8,6}$}}(6,-1)
\Label\o{\hbox{\scriptsize$X\!\!_{8,7}$}}(7,-1)
\Label\o{\hbox{\scriptsize$X\!\!_{8,8}$}}(8,-1)
\Label\o{\hbox{\scriptsize$X\!\!_{8,9}$}}(9,-1)
\Label\o{\hbox{\scriptsize$X\!\!_{8,1\!0}$}}(10,-1)
\Label\o{\hbox{\scriptsize$X\!\!_{8,1\!1}$}}(11,-1)
\Label\o{\hbox{\scriptsize$X\!\!_{8,1\!2}$}}(12,-1)
\Label\o{\hbox{\scriptsize$X\!\!_{8,1\!3}$}}(13,-1)
\Label\o{\hbox{\scriptsize$0$}}(0,0)
\Label\o{\hbox{\scriptsize$0$}}(1,0)
\Label\o{\hbox{\scriptsize$0$}}(2,0)
\Label\o{\hbox{\scriptsize$0$}}(3,0)
\Label\o{\hbox{\scriptsize$X\!\!_{7,4}$}}(4,0)
\Label\o{\hbox{\scriptsize$X\!\!_{7,5}$}}(5,0)
\Label\o{\hbox{\scriptsize$X\!\!_{7,6}$}}(6,0)
\Label\o{\hbox{\scriptsize$X\!\!_{7,7}$}}(7,0)
\Label\o{\hbox{\scriptsize$X\!\!_{7,8}$}}(8,0)
\Label\o{\hbox{\scriptsize$X\!\!_{7,9}$}}(9,0)
\Label\o{\hbox{\scriptsize$X\!\!_{7,1\!0}$}}(10,0)
\Label\o{\hbox{\scriptsize$X\!\!_{7,1\!1}$}}(11,0)
\Label\o{\hbox{\scriptsize$X\!\!_{7,1\!2}$}}(12,0)
\Label\o{\hbox{\scriptsize$X\!\!_{7,1\!3}$}}(13,0)
\Label\o{\hbox{\scriptsize$0$}}(0,1)
\Label\o{\hbox{\scriptsize$0$}}(1,1)
\Label\o{\hbox{\scriptsize$0$}}(2,1)
\Label\o{\hbox{\scriptsize$0$}}(3,1)
\Label\o{\hbox{\scriptsize$X\!\!_{6,4}$}}(4,1)
\Label\o{\hbox{\scriptsize$X\!\!_{6,5}$}}(5,1)
\Label\o{\hbox{\scriptsize$X\!\!_{6,6}$}}(6,1)
\Label\o{\hbox{\scriptsize$X\!\!_{6,7}$}}(7,1)
\Label\o{\hbox{\scriptsize$X\!\!_{6,8}$}}(8,1)
\Label\o{\hbox{\scriptsize$X\!\!_{6,9}$}}(9,1)
\Label\o{\hbox{\scriptsize$X\!\!_{6,1\!0}$}}(10,1)
\Label\o{\hbox{\scriptsize$X\!\!_{6,1\!1}$}}(11,1)
\Label\o{\hbox{\scriptsize$X\!\!_{6,1\!2}$}}(12,1)
\Label\o{\hbox{\scriptsize$X\!\!_{6,1\!3}$}}(13,1)
\Label\o{\hbox{\scriptsize$0$}}(0,2)
\Label\o{\hbox{\scriptsize$0$}}(1,2)
\Label\o{\hbox{\scriptsize$0$}}(2,2)
\Label\o{\hbox{\scriptsize$0$}}(3,2)
\Label\o{\hbox{\scriptsize$0$}}(4,2)
\Label\o{\hbox{\scriptsize$X\!\!_{5,5}$}}(5,2)
\Label\o{\hbox{\scriptsize$X\!\!_{5,6}$}}(6,2)
\Label\o{\hbox{\scriptsize$X\!\!_{5,7}$}}(7,2)
\Label\o{\hbox{\scriptsize$X\!\!_{5,8}$}}(8,2)
\Label\o{\hbox{\scriptsize$X\!\!_{5,9}$}}(9,2)
\Label\o{\hbox{\scriptsize$X\!\!_{5,1\!0}$}}(10,2)
\Label\o{\hbox{\scriptsize$X\!\!_{5,1\!1}$}}(11,2)
\Label\o{\hbox{\scriptsize$X\!\!_{5,1\!2}$}}(12,2)
\Label\o{\hbox{\scriptsize$X\!\!_{5,1\!3}$}}(13,2)
\Label\o{\hbox{\scriptsize$0$}}(0,3)
\Label\o{\hbox{\scriptsize$0$}}(1,3)
\Label\o{\hbox{\scriptsize$0$}}(2,3)
\Label\o{\hbox{\scriptsize$0$}}(3,3)
\Label\o{\hbox{\scriptsize$0$}}(4,3)
\Label\o{\hbox{\scriptsize$0$}}(5,3)
\Label\o{\hbox{\scriptsize$X\!\!_{4,6}$}}(6,3)
\Label\o{\hbox{\scriptsize$X\!\!_{4,7}$}}(7,3)
\Label\o{\hbox{\scriptsize$X\!\!_{4,8}$}}(8,3)
\Label\o{\hbox{\scriptsize$X\!\!_{4,9}$}}(9,3)
\Label\o{\hbox{\scriptsize$X\!\!_{4,1\!0}$}}(10,3)
\Label\o{\hbox{\scriptsize$X\!\!_{4,1\!1}$}}(11,3)
\Label\o{\hbox{\scriptsize$X\!\!_{4,1\!2}$}}(12,3)
\Label\o{\hbox{\scriptsize$X\!\!_{4,1\!3}$}}(13,3)
\Label\o{\hbox{\scriptsize$0$}}(0,4)
\Label\o{\hbox{\scriptsize$0$}}(1,4)
\Label\o{\hbox{\scriptsize$0$}}(2,4)
\Label\o{\hbox{\scriptsize$0$}}(3,4)
\Label\o{\hbox{\scriptsize$0$}}(4,4)
\Label\o{\hbox{\scriptsize$0$}}(5,4)
\Label\o{\hbox{\scriptsize$0$}}(6,4)
\Label\o{\hbox{\scriptsize$X\!\!_{3,7}$}}(7,4)
\Label\o{\hbox{\scriptsize$X\!\!_{3,8}$}}(8,4)
\Label\o{\hbox{\scriptsize$X\!\!_{3,9}$}}(9,4)
\Label\o{\hbox{\scriptsize$X\!\!_{3,1\!0}$}}(10,4)
\Label\o{\hbox{\scriptsize$X\!\!_{3,1\!1}$}}(11,4)
\Label\o{\hbox{\scriptsize$X\!\!_{3,1\!2}$}}(12,4)
\Label\o{\hbox{\scriptsize$X\!\!_{3,1\!3}$}}(13,4)
\Label\o{\hbox{\scriptsize$0$}}(0,5)
\Label\o{\hbox{\scriptsize$0$}}(1,5)
\Label\o{\hbox{\scriptsize$0$}}(2,5)
\Label\o{\hbox{\scriptsize$0$}}(3,5)
\Label\o{\hbox{\scriptsize$0$}}(4,5)
\Label\o{\hbox{\scriptsize$0$}}(5,5)
\Label\o{\hbox{\scriptsize$0$}}(6,5)
\Label\o{\hbox{\scriptsize$0$}}(7,5)
\Label\o{\hbox{\scriptsize$X\!\!_{2,8}$}}(8,5)
\Label\o{\hbox{\scriptsize$X\!\!_{2,9}$}}(9,5)
\Label\o{\hbox{\scriptsize$X\!\!_{2,1\!0}$}}(10,5)
\Label\o{\hbox{\scriptsize$X\!\!_{2,1\!1}$}}(11,5)
\Label\o{\hbox{\scriptsize$X\!\!_{2,1\!2}$}}(12,5)
\Label\o{\hbox{\scriptsize$X\!\!_{2,1\!3}$}}(13,5)
\Label\o{\hbox{\scriptsize$0$}}(0,6)
\Label\o{\hbox{\scriptsize$0$}}(1,6)
\Label\o{\hbox{\scriptsize$0$}}(2,6)
\Label\o{\hbox{\scriptsize$0$}}(3,6)
\Label\o{\hbox{\scriptsize$0$}}(4,6)
\Label\o{\hbox{\scriptsize$0$}}(5,6)
\Label\o{\hbox{\scriptsize$0$}}(6,6)
\Label\o{\hbox{\scriptsize$0$}}(7,6)
\Label\o{\hbox{\scriptsize$X\!\!_{1,8}$}}(8,6)
\Label\o{\hbox{\scriptsize$X\!\!_{1,9}$}}(9,6)
\Label\o{\hbox{\scriptsize$X\!\!_{1,1\!0}$}}(10,6)
\Label\o{\hbox{\scriptsize$X\!\!_{1,1\!1}$}}(11,6)
\Label\o{\hbox{\scriptsize$X\!\!_{1,1\!2}$}}(12,6)
\Label\o{\hbox{\scriptsize$X\!\!_{1,1\!3}$}}(13,6)
\Label\o{\hbox{\scriptsize$0$}}(0,7)
\Label\o{\hbox{\scriptsize$0$}}(1,7)
\Label\o{\hbox{\scriptsize$0$}}(2,7)
\Label\o{\hbox{\scriptsize$0$}}(3,7)
\Label\o{\hbox{\scriptsize$0$}}(4,7)
\Label\o{\hbox{\scriptsize$0$}}(5,7)
\Label\o{\hbox{\scriptsize$0$}}(6,7)
\Label\o{\hbox{\scriptsize$0$}}(7,7)
\Label\o{\hbox{\scriptsize$X\!\!_{0,8}$}}(8,7)
\Label\o{\hbox{\scriptsize$X\!\!_{0,9}$}}(9,7)
\Label\o{\hbox{\scriptsize$X\!\!_{0,1\!0}$}}(10,7)
\Label\o{\hbox{\scriptsize$X\!\!_{0,1\!1}$}}(11,7)
\Label\o{\hbox{\scriptsize$X\!\!_{0,1\!2}$}}(12,7)
\Label\o{\hbox{\scriptsize$X\!\!_{0,1\!3}$}}(13,7)
\hskip8.2cm
\right)
$$
\caption{}
\end{figure}

\begin{figure}[h] \label{fig2}
$$\Einheit.3cm
\Koordinatenachsen(14,16)(0,0)
\DuennPunkt(0,0)
\DuennPunkt(1,0)
\DuennPunkt(2,0)
\DuennPunkt(3,0)
\DuennPunkt(4,0)
\DuennPunkt(5,0)
\DuennPunkt(6,0)
\DuennPunkt(7,0)
\DuennPunkt(8,0)
\DuennPunkt(9,0)
\DuennPunkt(10,0)
\DuennPunkt(11,0)
\DuennPunkt(12,0)
\DuennPunkt(13,0)
\DuennPunkt(0,1)
\DuennPunkt(1,1)
\DuennPunkt(2,1)
\DuennPunkt(3,1)
\DuennPunkt(4,1)
\DuennPunkt(5,1)
\DuennPunkt(6,1)
\DuennPunkt(7,1)
\DuennPunkt(8,1)
\DuennPunkt(9,1)
\DuennPunkt(10,1)
\DuennPunkt(11,1)
\DuennPunkt(12,1)
\DuennPunkt(13,1)
\DuennPunkt(0,2)
\DuennPunkt(1,2)
\DuennPunkt(2,2)
\DuennPunkt(3,2)
\DuennPunkt(4,2)
\DuennPunkt(5,2)
\DuennPunkt(6,2)
\DuennPunkt(7,2)
\DuennPunkt(8,2)
\DuennPunkt(9,2)
\DuennPunkt(10,2)
\DuennPunkt(11,2)
\DuennPunkt(12,2)
\DuennPunkt(13,2)
\DuennPunkt(0,3)
\DuennPunkt(1,3)
\DuennPunkt(2,3)
\DuennPunkt(3,3)
\DuennPunkt(4,3)
\DuennPunkt(5,3)
\DuennPunkt(6,3)
\DuennPunkt(7,3)
\DuennPunkt(8,3)
\DuennPunkt(9,3)
\DuennPunkt(10,3)
\DuennPunkt(11,3)
\DuennPunkt(12,3)
\DuennPunkt(13,3)
\DuennPunkt(0,4)
\DuennPunkt(1,4)
\DuennPunkt(2,4)
\DuennPunkt(3,4)
\DuennPunkt(4,4)
\DuennPunkt(5,4)
\DuennPunkt(6,4)
\DuennPunkt(7,4)
\DuennPunkt(8,4)
\DuennPunkt(9,4)
\DuennPunkt(10,4)
\DuennPunkt(11,4)
\DuennPunkt(12,4)
\DuennPunkt(13,4)
\DuennPunkt(0,5)
\DuennPunkt(1,5)
\DuennPunkt(2,5)
\DuennPunkt(3,5)
\DuennPunkt(4,5)
\DuennPunkt(5,5)
\DuennPunkt(6,5)
\DuennPunkt(7,5)
\DuennPunkt(8,5)
\DuennPunkt(9,5)
\DuennPunkt(10,5)
\DuennPunkt(11,5)
\DuennPunkt(12,5)
\DuennPunkt(13,5)
\DuennPunkt(0,6)
\DuennPunkt(1,6)
\DuennPunkt(2,6)
\DuennPunkt(3,6)
\DuennPunkt(4,6)
\DuennPunkt(5,6)
\DuennPunkt(6,6)
\DuennPunkt(7,6)
\DuennPunkt(8,6)
\DuennPunkt(9,6)
\DuennPunkt(10,6)
\DuennPunkt(11,6)
\DuennPunkt(12,6)
\DuennPunkt(13,6)
\DuennPunkt(4,7)
\DuennPunkt(5,7)
\DuennPunkt(6,7)
\DuennPunkt(7,7)
\DuennPunkt(8,7)
\DuennPunkt(9,7)
\DuennPunkt(10,7)
\DuennPunkt(11,7)
\DuennPunkt(12,7)
\DuennPunkt(13,7)
\DuennPunkt(4,8)
\DuennPunkt(5,8)
\DuennPunkt(6,8)
\DuennPunkt(7,8)
\DuennPunkt(8,8)
\DuennPunkt(9,8)
\DuennPunkt(10,8)
\DuennPunkt(11,8)
\DuennPunkt(12,8)
\DuennPunkt(13,8)
\DuennPunkt(4,9)
\DuennPunkt(5,9)
\DuennPunkt(6,9)
\DuennPunkt(7,9)
\DuennPunkt(8,9)
\DuennPunkt(9,9)
\DuennPunkt(10,9)
\DuennPunkt(11,9)
\DuennPunkt(12,9)
\DuennPunkt(13,9)
\DuennPunkt(5,10)
\DuennPunkt(6,10)
\DuennPunkt(7,10)
\DuennPunkt(8,10)
\DuennPunkt(9,10)
\DuennPunkt(10,10)
\DuennPunkt(11,10)
\DuennPunkt(12,10)
\DuennPunkt(13,10)
\DuennPunkt(6,11)
\DuennPunkt(7,11)
\DuennPunkt(8,11)
\DuennPunkt(9,11)
\DuennPunkt(10,11)
\DuennPunkt(11,11)
\DuennPunkt(12,11)
\DuennPunkt(13,11)
\DuennPunkt(7,12)
\DuennPunkt(8,12)
\DuennPunkt(9,12)
\DuennPunkt(10,12)
\DuennPunkt(11,12)
\DuennPunkt(12,12)
\DuennPunkt(13,12)
\DuennPunkt(8,13)
\DuennPunkt(9,13)
\DuennPunkt(10,13)
\DuennPunkt(11,13)
\DuennPunkt(12,13)
\DuennPunkt(13,13)
\DuennPunkt(8,14)
\DuennPunkt(9,14)
\DuennPunkt(10,14)
\DuennPunkt(11,14)
\DuennPunkt(12,14)
\DuennPunkt(13,14)
\DuennPunkt(8,15)
\DuennPunkt(9,15)
\DuennPunkt(10,15)
\DuennPunkt(11,15)
\DuennPunkt(12,15)
\DuennPunkt(13,15)
\Label\lu{\scriptstyle0}(0,0)
\Label\u{\scriptstyle1}(1,0)
\Label\u{\scriptstyle2}(2,0)
\Label\u{\scriptstyle3}(3,0)
\Label\u{\scriptstyle4}(4,0)
\Label\u{\scriptstyle5}(5,0)
\Label\u{\scriptstyle6}(6,0)
\Label\u{\scriptstyle7}(7,0)
\Label\u{\scriptstyle8}(8,0)
\Label\u{\scriptstyle9}(9,0)
\Label\u{\scriptstyle10}(10,0)
\Label\u{\scriptstyle11}(11,0)
\Label\u{\scriptstyle12}(12,0)
\Label\u{\scriptstyle13}(13,0)
\Label\l{\scriptstyle1}(0,1)
\Label\l{\scriptstyle2}(0,2)
\Label\l{\scriptstyle3}(0,3)
\Label\l{\scriptstyle4}(0,4)
\Label\l{\scriptstyle5}(0,5)
\Label\l{\scriptstyle6}(0,6)
\Label\l{\scriptstyle7}(0,7)
\Label\l{\scriptstyle8}(0,8)
\Label\l{\scriptstyle9}(0,9)
\Label\l{\scriptstyle10}(0,10)
\Label\l{\scriptstyle11}(0,11)
\Label\l{\scriptstyle12}(0,12)
\Label\l{\scriptstyle13}(0,13)
\Label\l{\scriptstyle14}(0,14)
\Label\l{\scriptstyle15}(0,15)
\hskip4.2cm
$$
\caption{}
\end{figure}

Now fix a ``bivector" $M=[u_1,u_2,\dots,u_n\mid v_1,v_2,\dots,v_n]$ of
positive integers with $u_1<u_2<\dots<u_n$ and $v_1<v_2<\dots<v_n$.
Let $K[Y]$ denote the ring of all polynomials over some field $K$ in
the $Y_{i,j}$'s, $0\le i\le b$, $0\le j\le a$, and let $I_M(Y)$ be the
ideal in $K[Y]$ that is generated by those $t\times t$ minors of $Y$
that contain only nonzero entries, whose rows form a subset of the
last $u_t-1$ rows, or whose columns form a subset of the last $v_t-1$
columns, $t=1,2,\dots,n+1$. Here, by convention, $u_{n+1}$ is set
equal to $b+2$, and $v_{n+1}$ is set equal to $a+2$.  (Thus, for
$t=n+1$ the rows and columns of minors are unrestricted.)  The ideal
$I_M(Y)$ is called a {\it ladder determinantal ideal cogenerated by
the minor $M$}. (That one speaks of `the minor $M$' has its
explanation in the identification of the bivector $M$ with a
particular minor of $Y$, cf\@. \cite[Sec.~2]{HeTrAA}. It should be
pointed out that our conventions here deviate slightly from the ones
in \cite{HeTrAA}. In particular, we defined the ideal $I_M(Y)$ by
restricting rows and columns of minors to a certain number of {\it
last\/} rows or columns, while in \cite{HeTrAA} it is {\it first\/}
rows, respectively columns. Clearly, a rotation of the matrix by $180^\circ$
transforms one convention into the other.)  The associated {\it ladder
determinantal ring cogenerated by $M$\/} is $R_M(Y):=K[Y]/I_M(Y)$. (We
remark that the definition of ladder is more general in 
\cite{AbhyAB,AbKuAC,ConcAB,HeTrAA}. However, there is in effect no loss of
generality since the ladders of \cite{AbhyAB,AbKuAC,ConcAB,HeTrAA} 
can always be reduced to our definition by discarding
superfluous 0's.)

When Abhyankar introduced ladder determinantal rings in the
early\break 1980s, 
he was motivated by the study of singularities of Schubert varieties. 
Indeed, as was shown recently by
Gonciulea and Lakshmibai in \cite{GoLaAA} (see also 
\cite[Ch.~12]{BiLaAA}), the associated varieties (called ladder
determinantal varieties) can be identified with opposite cells of
certain Schubert varieties of type $A$. This connection allowed them
to identify the irreducible components of such Schubert
varieties in many cases, thus making substantial progress on a
long-standing problem in algebraic geometry.

Results of Abhyankar 
\cite{AbhyAB,AbKuAC} or Herzog and Trung \cite{HeTrAA}
allow to express the
Hilbert series of the ladder determinantal ring $R_M(Y)$ in 
combinatorial terms. Before we can state the corresponding result, we
need to introduce a few more terms.

When we say {\it lattice path\/} we always mean
a lattice path in the plane consisting of unit horizontal and vertical
steps in the positive direction, see Figure~3 for an example. 
We shall frequently abbreviate the fact that a lattice path $P$ goes from
$A$ to $E$ by $P:A\to E$.

\begin{figure}[h] \label{fig3}
$$\Koordinatenachsen(8,8)(0,0)
\Gitter(8,8)(0,0)
\Pfad(1,-1),221221112122\endPfad
\DickPunkt(1,-1)
\DickPunkt(6,6)
\Label\ro{P_0}(3,3)
\hskip4cm
$$
\caption{}
\end{figure}

Also, given
lattice points $A$ and $E$, we denote the set of all lattice paths
from $A$ to $E$ by $\P(A\to E)$.  A family $(P_1,P_2,\dots,P_n)$ of
lattice paths $P_{i}$, $i=1,2,\dots,n$, is said to be {\it
nonintersecting\/} if no two lattice paths of this family have a point
in common.  Given $n$-tuples of lattice points $\mathbf
A=(A^{(1)},A^{(2)},\dots,A^{(n)})$ and 
$\mathbf E=(E^{(1)},E^{(2)},\dots,E^{(n)})$, we denote
the set of all families $(P_1,P_2,\dots,P_n)$ of nonintersecting
lattice paths, where $P_i$ runs from $A^{(i)}$ to $E^{(i)}$,
$i=1,2,\dots,n$, by $\P^+(\mathbf A\to\mathbf E)$.

A point in a lattice path $P$ which is the end point of a vertical
step and at the same time the starting point of a horizontal step will
be called a {\it north-east turn\/} ({\it NE-turn\/} for short) of the
lattice path $P$.  The NE-turns of the lattice path in
Figure~3 are
$(1,1)$, $(2,3)$, and $(5,4)$. We write $\NE(P)$ for the number of
NE-turns of $P$. Also, given a family $\mathbf P=(P_1,P_2,\dots,P_n)$ of
paths $P_i$, we write $\NE(\mathbf P)$ for the number $\sum _{i=1}
^{n}\NE(P_i)$ of all NE-turns in the family.

Our lattice paths will be restricted to ladder-shaped regions $L$
corresponding to the nonzero entries of a given matrix $Y$ in the way
that was explained earlier (cf\@. Figures~1 and 2).  
We extend our
lattice path notation in the following way. By $\P_L(A\to E)$ we mean
the set of all lattice paths $P$ from $A$ to $E$ {\it all of whose
NE-turns lie in the ladder region $L$}.  (It should be noted that, in
the case of a two-sided ladder, it is possible that a path is {\it
not\/} totally inside $L$ while its NE-turns {\it are}. However, in
the case of an upper ladder $L$, which is the case of interest for our
main results Theorem~\ref{thm:2} and Corollary~\ref{cor:3}, a path is inside $L$ {\it if
and only if\/} all of its NE-turns are.)  Similarly, by $\P^+_L(\mathbf
A\to \mathbf E)$ we mean the set of all families $(P_1,P_2,\dots,P_n)$
of nonintersecting lattice paths, where $P_i$ runs from $A^{(i)}$ to
$E^{(i)}$ {\it and where all the NE-turns of $P_i$ lie in the ladder
region $L$}.

Finally, given any weight function $w$ defined on a set $\mathcal M$, by
the generating function $\GF(\mathcal M;w)$ we mean $\sum _{x\in\mathcal M}
^{}w(x)$.

\begin{theorem}
\label{thm:1}
Let $Y=(Y_{i,j})_{0\le i\le b,\ 0\le j\le a}$ be
a (two-sided) ladder, and let $L$ be the associated ladder region, 
i.e., $Y_{i,j}=X_{i,j}$ if and only
if $(j,b-i)\in L$. Let $M=[u_1,u_2,\dots,u_n\mid v_1,v_2,\dots,v_n]$
be a bivector of positive integers with 
$u_1<u_2<\dots<u_n$ and $v_1<v_2<\dots<v_n$.
Furthermore, let $A^{(i)}=(0,u_{n+1-i}-1)$ and $E^{(i)}=(a-v_{n+1-i}+1,b)$, 
$i=1,2,\dots,n$. 
Then, under the assumption that all of the points $A^{(i)}$ and
$E^{(i)}$, $i=1,2,\dots,n$,
lie inside the ladder region $L$, the Hilbert series of the ladder determinantal
ring $R_M(Y)=K[Y]/I_M(Y)$ equals
\begin{equation}
\sum _{\ell=0} ^{\infty}\dim_K R_M(Y)_\ell\,z^\ell
=\frac {\GF(\P_L^+(\mathbf A\to \mathbf E);z^{\NE(.)})}
{(1-z)^{(a+b+3)n-\sum _{i=1} ^{n}(u_i+v_i)}},
\label{e2.1}
\end{equation}
where $R_M(Y)_\ell$ denotes the homogeneous component of degree $\ell$
in $R_M(Y)$, and where, according to our definitions,
$\GF(\P_L^+(\mathbf A\to \mathbf E);z^{\NE(.)})$ is the generating
function $\sum _{\mathbf P} ^{}z^{\NE(\mathbf P)}$ for all families $\mathbf
P=(P_1,P_2,\dots,P_n)$ of nonintersecting lattice paths, $P_i$ running
from $A^{(i)}$ to $E^{(i)}$, such that all of its NE-turns stay inside
the ladder region $L$.
\end{theorem}

\begin{remark} \label{rem:1}
The condition that all of the points $A^{(i)}$ and $E^{(i)}$
lie inside the ladder region $L$ restricts the choice of ladders. In
particular, for an upper ladder it means that 
$Y_{b-u_{n}+1,0}=X_{b-u_{n}+1,0}$ and
$Y_{0,a-v_{n}+1}=X_{0,a-v_{n}+1}$, which will be relevant for us.
Still, one could prove an analogous result
even if this condition is dropped. In that case, however,
the points $A^{(i)}$ and $E^{(i)}$ have to be modified in
order to lie inside $L$ and, thus, make the right-hand side of
formula \eqref{e2.1} meaningful. 
\end{remark}

\medskip
\noindent
{\it Sketch of Proof}.
In \cite[proof of Theorem~2]{KrPrAA}, we gave
two proofs of this assertion in the special case of a one-sided
ladder and $u_i=v_i=i$, $i=1,2,\dots,n$ (cf\@. Example~(1) on p.~10 of
\cite{HeTrAA}). The first proof followed basically considerations by
Kulkarni \cite{KulkAD,KulkAC} (see also \cite{GhorAD}), and was
based on an explicit basis for $R_M(Y)$ given by Abhyankar
\cite[Theorem~(20.10)(5)]{AbhyAB}.  The second proof was based on
combinatorial descriptions of the dimensions $R_M(Y)_\ell$ of the
homogeneous components of $R_M(Y)$ due to Herzog and Trung
\cite[Cor.~4.3 + Lemma~4.4]{HeTrAA}.  Both proofs carry over verbatim
to our more general situation because both Abhyankar's as well as
Herzog and Trung's results are in fact theorems for the general ladder
determinantal rings that we consider here. (However, the reader must
be warned that the explicit form of Abhyankar's basis was misquoted in
\cite{KrPrAA}. The correct assertion is that, given a multiset $S$ as
described in \cite{KrPrAA}, the associated basis element is the
product of a certain monomial in the $X_{ij}$'s 
and a certain minor of the matrix $Y$,
see \cite[definition of $w_v(t)$ in Theorem~(20.10)]{AbhyAB} 
or \cite[Theorem~(6.7)(iii)]{GhorAD}
Also, the definition of the multisets $S$ contained an error: Item 2
at the bottom of p.~1019 in \cite{KrPrAA} must be replaced by: The
length of any sequence $(i_1,j_1)$, $(i_2,j_2)$, \dots, $(i_k,j_k)$
of elements of $S$ is at most $n$. The subsequent argument was
however based on this corrected definition.) 
\quad \quad \qed

\section{The determinantal formula}
\label{Sec:3}

In view of Theorem~\ref{thm:1}, the computation of Hilbert series of ladder
determinantal rings requires to solve the problem of counting families
of nonintersecting lattice paths in a ladder-shaped region with
respect to turns. We provide such a solution for one-sided ladders in
Theorem~\ref{thm:2}. In order to formulate the result, we need to introduce the
notion of {\it two-rowed arrays}.

{}From now on we restrict our attention to one-sided ladders.
Without loss of generality it suffices to consider upper ladders.
We encode upper ladder-shaped regions (such as the one in
Figure~2) 
concisely by means of weakly
increasing functions as follows: 
given an upper ladder region $L$,
let $f$ be the weakly increasing function from $[0,a]$ to
$[1,b+1]$ with the property that it describes $L$ by means of
\begin{equation}
L=\{(x,y):x\in[0,a]\text{ and }0\le y<f(x)\}.\label{e3.1}
\end{equation} 
Here, by
$[c,d]$ we mean the set of all integers $\ge c$ and $\le d$. In
essence, the function $f$ describes the upper border of the region
$L$. For example, the function $f$ corresponding to the ladder region
in Figure~2 (where $a=13$ and $b=15$) is given by
$f(0)=7$, $f(1)=7$, $f(2)=7$, $f(3)=7$, $f(4)=10$, $f(5)=11$,
$f(6)=12$, 
$f(7)=13$, $f(8)=16$, $f(9)=16$, $f(10)=16$, $f(11)=16$, $f(12)=16$,
$f(13)=16$. 

By a {\it two-rowed array} we mean two rows of integers
\begin{equation}
\begin{matrix}
 a_{-l+1}&a_{-l+2}&\dots&a_{-1}&a_0&a_1&\dots&a_k\phantom{,}\\
         &        &     &      &   &b_1&\dots&b_k,
\end{matrix}\label{e3.2}
\end{equation}
where entries along both rows are {\it strictly increasing}. We call
$l$ the {\it type} of the two-rowed array. We allow $l$ to
be also negative. In this case the representation \eqref{e3.2} has to be taken
symbolically, in the sense that the first row of the two-rowed array
is (by $-l$) shorter than the second row, i.e., looks like
\begin{equation}\label{e3.2a}
\begin{matrix}
         &        &     &      &a_{-l+1}&\dots&a_k\phantom{,}\\
 b_1     &b_2     &\dots&b_{-l}&b_{-l+1}&\dots&b_k.
\end{matrix}
\end{equation}

We define the {\it size} $\vert T\vert$ of a two-rowed array $T$ to be
the number of its entries. (Thus, the size of the two-rowed array in
\eqref{e3.2} is $l+2k$, as is the size of the one in \eqref{e3.2a}.) 
We extend this definition and notation to families
$\mathbf T=(T_1,T_2,\dots,T_n)$ of two-arrays by letting $\vert\mathbf
T\vert$ denote the total number $\vert T_1\vert+\vert
T_2\vert+\dots+\vert T_n\vert$ of entries in $\mathbf T$.

Now we define the basic set of objects which is crucial in our
formulas.
Given a function $f$ as above, and pairs $A=(\al_1,\al_2)$ and
$E=(\ep_1,\ep_2)$, we denote by $\TA(l;A,E;f,d)$ the set of all
two-rowed arrays of type $l$ such that
\begin{itemize}
\item the entries in the first row are bounded below by
$\al_1$ and bounded above by $\ep_1$,
\item the entries in the second row are bounded below by
$\al_2$ and bounded above by $\ep_2$,
\item if the two-rowed array is represented as in \eqref{e3.2}
(respectively \eqref{e3.2a}), we
have
\begin{equation}
b_s<f(a_{s+d}), \label{e3.3}
\end{equation}
for all $s$ such that both $b_s$ and $a_{s+d}$ exist in the two-rowed
array. 
\end{itemize}

If we want to make the lower and upper bounds transparent,
then we will write such two-rowed arrays in the form
\begin{equation}
\begin{array}{rccccccccl}
\al_1\le&a_{-l+1}&a_{-l+2}&\dots&a_{-1}&a_0&a_1&\dots&a_k&\le\ep_1\\
\al_2\le&        &        &     &      &   &b_1&\dots&b_k&\le\ep_2.
\end{array}\label{e3.4}
\end{equation}

Our key theorem is the following.  
\begin{theorem}
\label{thm:2}
Let $n,a,b$ be positive integers and let $L$ be an upper ladder-shaped
region determined by the weakly increasing function
$f:[0,a]\to[1,b+1]$ by means of {\em \eqref{e3.1}}. For convenience,
extend $f$ to all negative integers by setting $f(x):=f(0)$ for
$x<0$. Furthermore, let
$A^{(i)}=(A^{(i)}_1,A^{(i)}_2)$ and $E^{(i)}=(E^{(i)}_1,E^{(i)}_2)$ for
$i=1,2,\dots,n$ be lattice points in the region $L$ satisfying 
\begin{equation} \label{e3.5c}
f(x)=f\big(A^{(1)}_1\big)\quad \text {for all $x\le A^{(1)}_1$,}
\end{equation}
and
\begin{gather} \label{e3.5a}
A^{(1)}_1\le A^{(2)}_1\le\dots\le A^{(n)}_1,\quad
A^{(1)}_2 > A^{(2)}_2 > \dots > A^{(n)}_2,
\end{gather}
and
\begin{gather} \label{e3.5b}
E^{(1)}_1 < E^{(2)}_1 < \dots < E^{(n)}_1,\quad
E^{(1)}_2\ge E^{(2)}_2\ge\dots\ge E^{(n)}_2.
\end{gather}
Then the generating function $\sum z^{\NE(\mathbf P)}$, where the sum is
over all families $\mathbf P =(P_1,P_2,\dots,P_n)$ of nonintersecting
lattice paths $P_i:A^{(i)}\to E^{(i)}$, $i=1,2,\dots,n$ lying in the
region $L$, can be expressed as
\begin{multline}
\GF(\P_L^+(\mathbf A\to\mathbf E);z^{\NE(.)})\\
=\det_{1\le s,t\le n}
   \big(\GF(\TA(t-s;\tilde A^{(t)},\tilde E^{(s)};f,s-1);z^{\vert.\vert/2})\big),
\label{e3.5}
\end{multline}
where $\tilde A^{(i)}=A^{(i)}+(-i+1,i)$ and $\tilde
E^{(i)}=E^{(i)}+(-i,i-1)$, $i=1,2,\dots,n$. Here, by our definitions,
$\GF(\TA(t-s;\tilde A^{(t)},\tilde E^{(s)};f,s-1);z^{\vert.\vert/2})$ is the
generating function $\sum_{T}z^{\vert T\vert/2}$, where the sum is
over all two-rowed arrays of the form {\em \eqref{e3.4}} with $l=t-s$,
$d=s-1$, $\al_1=A^{(i)}_1-i+1$, $\al_2=A^{(i)}_2+i$, $\ep_1=E^{(i)}_1-i$,
and $\ep_2=E^{(i)}_2+i-1$, which satisfy {\em \eqref{e3.3}}.
\end{theorem}

\begin{remark} (1) The condition \eqref{e3.5c} is equivalent to saying 
that to the left of
$A^{(1)}$, which by \eqref{e3.5a} is
the left-most starting point of the lattice paths, 
the boundary of the ladder region is horizontal. Clearly, this can be
assumed without loss of generality because this part of the ladder
(i.e., the ladder to the left of $A^{(1)}$) does not impose any
restriction on the lattice paths, and, hence, on the left-hand side
of \eqref{e3.5}.

(2) The formula \eqref{e3.5} clearly reduces the problem of enumerating 
families of nonintersecting lattice paths in the ladder region $L$
with respect to NE-turns to the problem of enumerating certain
two-rowed arrays. We are going to address 
this problem in Section~\ref{Sec:5}.
\end{remark}

Thus, if we combine Theorems~\ref{thm:1} and \ref{thm:2}, we obtain the promised
determinantal formula for the Hilbert series of one-sided ladder
determinantal rings.

\begin{corollary}
\label{cor:3}
Let $Y=(Y_{i,j})_{0\le i\le b,\ 0\le j\le a}$ be an upper ladder, and
let $L$ be the associated ladder region, i.e., $Y_{i,j}=X_{i,j}$ if
and only if $(j,b-i)\in L$, and let $f:[0,a]\to [1,b+1]$ be the
function that describes this ladder region by means of {\em \eqref{e3.1}},
i.e., $Y_{i,j}=X_{i,j}$ if and only if $b-i<f(j)$. For convenience,
extend $f$ to all negative integers by setting $f(x):=f(0)$ for
$x<0$. Let
$M=[u_1,u_2,\dots,u_n\mid v_1,v_2,\dots,v_n]$ be a bivector of
positive integers with $u_1<u_2<\dots<u_n$ and $v_1<v_2<\dots<v_n$ such that
$Y_{b-u_{n}+1,0}=X_{b-u_{n}+1,0}$ and
$Y_{0,a-v_{n}+1}=X_{0,a-v_{n}+1}$ (cf\@. Remark~\ref{rem:1} after
Theorem~\ref{thm:1}).
Furthermore, we let $\tilde
A^{(i)}=(-i+1,u_{n+1-i}+i-1)$ and $\tilde
E^{(i)}=(a-v_{n+1-i}-i+1,b+i-1)$, $i=1,2,\dots,n$.  Then the Hilbert
series of the ladder determinantal ring $R_M(Y)=K[Y]/I_M(Y)$ equals
\begin{multline}
\sum _{\ell=0} ^{\infty}\dim_K R_M(Y)_\ell\,z^\ell\\
=\frac {\det_{1\le s,t\le n}
\big(\GF(\TA(t-s;\tilde A^{(t)},\tilde E^{(s)};f,s-1);
z^{\vert.\vert/2})\big)}{(1-z)^{(a+b+3)n-\sum_{i=1}^{n}(u_i+v_i)}}.
\label{e3.6}
\end{multline}
\end{corollary}

\begin{remark} (1) Theorem~\ref{thm:2} specializes to Theorem~1 in
\cite{KratBE} in the case of a trivial ladder (i.e., if
the function $f$ is equal to $b+1$ for all $x$). For, in that case, 
by \eqref{e5.4a}
the generating functions $\GF(\TA(t-s;\tilde A_t,\tilde E_s;f,s-1);
z^{\vert.\vert/2})$ can be expressed in terms of binomial sums.
To see that the resulting formula is indeed equivalent, one extracts
the coefficient of $z^K$.

(2) For the same reason, Corollary~\ref{cor:3} specializes to Abhyankar's
formula \cite[(20.14.4), $L=2$, $k=2$, with
$F^{(22)}(m,p,a,V)$ defined on p.~50]{AbhyAB} 
in the case of a trivial ladder. 
Although Abhyankar's formula gives an expression for the Hilbert
function (instead of for the Hilbert series), it is easy to see that
it is equivalent to ours in this special case.

(3) The formula for the Hilbert series in \cite[Theorem~2]{KrPrAA}
addresses the special case $u_i=v_i=i$, $i=1,2,\dots,n$. However,
Corollary~\ref{cor:3} does not generalize this formula, as it does
not directly specialize to Theorem~2 in \cite{KrPrAA}. Whereas in the
latter formula the entries of the determinant are generating
functions for paths, there is no such interpretation for the entries
of the determinant in \eqref{e3.6}.

(4) Unfortunately, we do not know how to generalize
Theorem~\ref{thm:2}, and, thus, Corollary~\ref{cor:3}, to the case of
two-sided ladders. It seems that a completely new idea is needed to
find such a generalization.

(5) More modest, but equally desirable, would it be to find an extension of 
Corollary~\ref{cor:3} in the
one-sided case to ladders $L$ and bivectors $M$ which do not satisfy
the conditions of the statement, i.e., for which either 
$Y_{b-u_{n}+1,0}=0$, or $Y_{0,a-v_{n}+1}=0$, or both.
This would require to find an extension of Theorem~\ref{thm:2} to
situations where the inequality chains \eqref{e3.5a} and \eqref{e3.5b}
may be relaxed so that some starting and end points are allowed to lie on
the boundary of the ladder region $L$ (cf\@. Remark~\ref{rem:1} after
Theorem~\ref{thm:1}).
It seems again that a completely new idea is needed to
find such an extension.

(6) In Section~5 of \cite{KrPrAA} it is shown that the proof of the
main counting theorem yields in fact a
weighted generalization thereof.
An analogous weighted generalization of Theorem~\ref{thm:2} can be
obtained as well, which is again directly implied by the proof of
Theorem~\ref{thm:2} in
Section~\ref{Sec:4}. However, we omit the statement of this
generalization for the sake of brevity.
\end{remark}

\begin{example} Let $a=13$, $b=15$, $n=4$, let
$Y=(Y_{i,j})$ be the matrix of Figure~1 and
$M=[1,2,4,6\mid 1,2,3,6]$. Our
formula \eqref{e3.6} gives for the Hilbert series of
$R_{M}(Y)=K[Y]/I_{M}(Y)$, using \eqref{e5.11} for determining the generating
function $\sum_{T\in\TA(l;A,E;f,d)}z^{\vert.\vert/2}$ for two-rowed
arrays $T$ in the corresponding ladder region $L$ of
Figure~2,
\begin{align*} \scriptstyle 
(1 &\scriptstyle+ 71 z + 2556 z^2 + 61832 z^3 + 1115762 z^4 + 15750005 z^5 +
    178390279 z^6 + 1647137174 z^7 \\
&\scriptstyle+ 12534233703 z^8+ 79245271879 z^9 + 418852424787 z^{10} + 1859941402206 z^{11}
    + 6965987806143 z^{12} \\
&\scriptstyle + 22071622313567 z^{13}+ 59298706514083 z^{14} + 135299444287353 z^{15} +
    262400571075662 z^{16} \\
&\scriptstyle + 432640455645309 z^{17}+ 606103694379729 z^{18} + 720535170430557 z^{19} +
    725289798304502 z^{20} \\
&\scriptstyle+ 616230022969392 z^{21}+ 439998448014899 z^{22} + 262469031030333 z^{23} +
    129776697745621 z^{24} \\
&\scriptstyle + 52622863698472 z^{25}+ 17241967478923 z^{26} + 4468021840695 z^{27} + 885721405230
    z^{28}  \\
&\scriptstyle + 126901720400 z^{29}+ 11760999250 z^{30}+ 532021875  z^{31} )/ (1-z)^{99}.
\end{align*}

\end{example}

\section{Proof of Theorem~\ref{thm:2}}
\label{Sec:4}

The basic idea of the proof is simple. It largely follows the proof
of Theorem~4 in \cite{KratBE}. As a first step, we expand
the determinant on the right-hand side of \eqref{e3.5} according to
the definition of a determinant, see Subsection~\ref{Sub:Exp}. Thus,
we obtain a sum of terms, each of which is indexed by a family of
two-rowed arrays, see \eqref{e4.2}. Some of the terms have positive
sign, some of them negative sign. In the second step, we identify the
terms which cancel each other, see Subsection~\ref{Sub:Cancel}. 
Finally, in the third step, we identify the remaining terms with the
families of nonintersecting lattice paths in the statement of the
theorem, see Subsection~\ref{Sub:Remain}. 

However, the details are sometimes intricate. To show that the terms
described in Subsection~\ref{Sub:Cancel} do indeed cancel, we define
an involution on families of two-rowed arrays in
Subsection~\ref{Sub:Inv}. (This involution is copied from
\cite[Proof of Theorem~4]{KratBE}.) In order that our claims follow,
this involution must have several properties, which are listed in
Subsection~\ref{Sub:Prop}. While most of these are either obvious or are
already established in \cite{KratBE} and \cite{RubeAB}, we are only
able to provide a rather technical
justification of the one pertaining to the ladder condition. This is
done in Subsection~\ref{Sub:Ladder}.

\subsection{Expansion of the determinant} \label{Sub:Exp}
Let $\S_n$ denote the symmetric group of order $n$.  We start by
expanding the determinant on the right-hand side of \eqref{e3.5}, to obtain
\begin{align} \notag
\det_{1\le s,t\le n}\big(\GF(&\TA(t-s;\tilde A^{(t)},\tilde E^{(s)};f,s-1);
z^{\vert.\vert/2})\big)\\
\notag
&=\sum_{\si\in \S_n}\sgn\si\prod_{i=1}^{n}
\GF(\TA(\si(i)-i;\tilde A^{(\si(i))},\tilde E^{(i)};f,i-1);z^{\vert.\vert/2})\\
&=\sum_{({\mathbf T},\si)}\sgn\si\,z^{\vert {\mathbf T}\vert},
\label{e4.2}
\end{align}
where the sum is over all pairs $({\mathbf T}, \si)$ of
permutations $\si$ in $\S_n$,
and families ${\mathbf T}=( T_1, T_2,\dots,
T_n)$ of two-rowed arrays, $ T_i$ being of type $\si(i)-i$ (i.e., the
second row containing $k_i$ entries and the first row containing
$k_i+\si(i)-i$ entries, for some $k_i$), and
the bounds for the entries of $ T_i$ being as follows,
\begin{equation}
\begin{array}{rccccl}
\tilde A^{(\si(i))}_1\le&a_{-\si(i)+i+1}^{(i)}\dots\ &a_1^{(i)}&\dots\ &a_{k_i}^{(i)}&\le \tilde E^{(i)}_1\\
\tilde A^{(\si(i))}_2\le&&b_1^{(i)}&\dots\ &b_{k_i}^{(i)}&\le \tilde E^{(i)}_2,
\end{array}\label{e4.3}
\end{equation}
with the property that 
\begin{equation}
b_s^{(i)}<f(a_{s+i-1}^{(i)}),\quad \quad s=1,2,\dots,k_i-i+1,\label{e4.4}
\end{equation}
$i=1,\dots,n$. 

\subsection{Which terms in \eqref{e4.2} cancel?}    \label{Sub:Cancel}
Now we claim that the total contribution to the sum \eqref{e4.2} of the
families $( T_1, T_2,\dots, T_n)$ of two-rowed
arrays as above which have the property that there exist $ T_i$
and $ T_{i+1}$, $ T_i$ represented by
{\refstepcounter{equation}\label{e4.5}}
\alphaeqn
\begin{equation}
\begin{array} {rccccl}
\tilde A^{(\si(i))}_1\le&a_{-\si(i)+i+1}\ \dots\ &a_1&\ \dots\ &a_{k}&\le \tilde E^{(i)}_1\\
\tilde A^{(\si(i))}_2\le&&b_1&\ \dots\ &b_{k}&\le \tilde E^{(i)}_2,
\end{array}\label{e4.5a}
\end{equation}
and $ T_{i+1}$ represented by
\begin{equation}
\begin{array} {rccccl}
\tilde A^{(\si(i+1))}_1\le&c_{-\si(i+1)+i+2}\ \dots\ &c_1&\ \dots\ &c_{l}&\le \tilde E^{(i+1)}_1\\
\tilde A^{(\si(i+1))}_1\le&&d_1&\ \dots\ &d_{l}&\le \tilde E^{(i+1)}_1,
\end{array}\label{e4.5b}
\end{equation}
and indices $I$ and $J$ such that
\begin{gather}  c_J<a_I\label{e4.5c}\\
b_{I-1}<d_J\label{e4.5d}
\end{gather}
and
\begin{equation}
1\le I\le k+1,\quad 0\le J\le l,\label{e4.5e}
\end{equation} 
\reseteqn
equals 0.  The
inequalities \eqref{e4.5c} and \eqref{e4.5d} should be understood to hold only if all
variables are defined, including the conventional definitions
$a_{k+1}:=\tilde E^{(i)}_1+1$, $b_0:=\tilde A^{(\si(i))}_2-1$, and
$c_{-\si(i+1)+i+1}:=\tilde A^{(\si(i))}_1-1$. (These artificial
settings apply for $I=k+1$, $I=1$, and $J=-\si(i+1)+i+1$,
respectively. It should be noted that the indexing conventions that we
have chosen here deviate slightly from \cite[Sec.~3, proof of
Theorem~4]{KratBE}, but are completely equivalent.)

We call the point $(a_I,d_J)$ a {\it crossing point\/} of $ T_i$
and $ T_{i+1}$, and, more generally, a {\it crossing point\/} of
the family ${\mathbf T}$.  

\subsection{The remaining terms correspond to nonintersecting lattice
paths} \label{Sub:Remain}
Suppose that we would have shown that the contribution to \eqref{e4.2} of
these families of two-rowed arrays equals zero. It implies that only
those families ${\mathbf T}=( T_1, T_2,\dots,
T_n)$ of two-rowed arrays, $ T_i$ being of the form \eqref{e4.3} and
satisfying \eqref{e4.4}, contribute to \eqref{e4.2} where $ T_i$ and $
T_{i+1}$ have no crossing point for all $i$.

So, let $\mathbf T$ be such a family of two-rowed arrays without any
crossing point. By using the arguments from 
\cite{RubeAB}\footnote{The proof
in the original paper \cite[last paragraph of the proof of
Theorem~4]{KratBE} contained an error at this point. The inequality
$A^{(\si(i+1))}_1-1\le A^{(\si(i))}$ on page~12 of \cite{KratBE} is
not true in general.} 
(with $A^{(i)}_1$, $A^{(i)}_2$, $E^{(i)}_1$, $E^{(i)}_2$
in \cite{RubeAB} replaced by our $\tilde A^{(i)}_1$, $\tilde A^{(i)}_2-1$,
$\tilde E^{(i)}_1+1$, $\tilde E^{(i)}_2$, respectively,
$i=1,2,\dots,n$),
it then follows that the permutation $\si$ associated
to ${\mathbf T}$ must be the identity permutation. 
Thus, the two-rowed array $ T_i$ has the form (recall \eqref{e4.3})
\begin{equation}
\begin{array} {rcccl}
\tilde A^{(i)}_1\le&a_1^{(i)}&\dots\ &a_{k_i}^{(i)}&\le \tilde E^{(i)}_1\\
\tilde A^{(i)}_2\le&b_1^{(i)}&\dots\ &b_{k_i}^{(i)}&\le \tilde E^{(i)}_2,
\end{array}\label{e4.6}
\end{equation}
and satisfies \eqref{e4.4}. Moreover, we assumed that
there is no crossing point, meaning
that there are no consecutive two-rowed arrays $ T_i$ and
$ T_{i+1}$ and indices $I$ and $J$ such that \eqref{e4.5} holds.

By interpreting the two-rowed array \eqref{e4.6} as a lattice path $\tilde
P_i$ from $\tilde A^{(i)}-(0,1)$ to $\tilde E^{(i)}+(1,0)$ whose
NE-turns are exactly $(a_1^{(i)},b_1^{(i)})$, \dots,
$(a_{k_i}^{(i)},b_{k_i}^{(i)})$, $i=1,2,\dots,n$, the family
${\mathbf T}$ of two-rowed arrays is translated into a family
$\widetilde{\mathbf P}=(\tilde P_1,\tilde P_2,\dots,\tilde P_n)$ of
paths. Clearly, under this translation we have $\vert{\mathbf
T}\vert/2=\NE(\widetilde{\mathbf P})$, and, hence,
\begin{equation}
z^{\vert{\mathbf T}\vert/2}=z^{\NE(\widetilde{\mathbf P})}.\label{e4.7}
\end{equation}

The fact that \eqref{e4.5} does not hold simply means that the paths $\tilde
P_i$ and $\tilde P_{i+1}$ do not cross each other (that is, they may
touch each other, but they never change sides), $i=1,2,\dots,n-1$. We
refer the reader to the explanations in Section~\ref{Sec:2} (between Theorems~3
and 4) in \cite{KratBE}. Here, we content ourselves with an
illustration.  Suppose two paths $Q_1$ and $Q_2$ cross each other (see
Figure~4).  Furthermore suppose that the NE-turns of $Q_1$ are
$(a_1,b_1)$, $(a_2,b_2)$, \dots, $(a_k,b_k)$, and the NE-turns of
$Q_2$ are $(c_1,d_1)$, $(c_2,d_2)$, \dots, $(c_l,d_l)$.  Then it is
obvious from Figure~4 that there exist $I$ and $J$ such that
\eqref{e4.5c}--\eqref{e4.5e} hold.

\begin{figure}[h] \label{fig4}
$$
\Pfad(0,0),221111121111\endPfad
\Pfad(-2,-1),11112222221111\endPfad
\DickPunkt(0,2)
\DickPunkt(2,2)
\DickPunkt(2,-1)
\Label\lo{\hskip-10pt(c_J,d_J)}(0,2)
\Label\ru{(a_I,b_{I-1})}(2,-1)
\Label\r{Q_1}(6,5)
\Label\r{Q_2}(9,3)
\hskip4cm
$$
\caption{}
\end{figure}

To finally match with the claim of Theorem~\ref{thm:2}, we shift $\tilde P_i$ by
$(i-1,-i+1)$, $i=1,2,\dots,n$. Thus we obtain a family
$(P_1,P_2,\dots,P_n)$ of lattice paths, $P_i$ running from $A^{(i)}$
to $E^{(i)}$. Clearly, under this shift, the condition that $\tilde
P_i$ and $\tilde P_{i+1}$ do not cross each other translates into the
condition that $P_i$ and $P_{i+1}$ do not touch each other,
$i=1,2,\dots,n-1$.  If we combine this fact with the observation that
the first path, $P_1=\tilde P_1$, stays inside the ladder region $L$
because of \eqref{e4.4} with $i=1$, then we conclude that all the $P_i$'s must also stay
inside $L$ because $P_1$ forms a barrier.

Thus, in view of \eqref{e4.7}, we have proved that the right-hand side of
\eqref{e3.5} is equal to the generating function $\sum_{\mathbf P} z^{\NE(\mathbf
P)}$, where the sum is over all families $\mathbf P=(P_1,P_2,\dots,P_n)$ of
nonintersecting lattice paths, $P_i$ running from $A^{(i)}$ to
$E^{(i)}$ and staying inside the ladder region $L$. But this is
exactly the left-hand side of \eqref{e3.5}.  Thus Theorem~\ref{thm:2} would be proved.

\subsection{The involution} \label{Sub:Inv}
To show that the contribution to the sum \eqref{e4.2} of the families
${\mathbf T}=( T_1, T_2,\dots,\break T_n)$ of two-rowed
arrays, $ T_i$ being of the form \eqref{e4.3} and satisfying \eqref{e4.4} for
$i=1,2,\dots,n$, which contain consecutive arrays $ T_i$ and
$ T_{i+1}$ that have a crossing point (cf\@. \eqref{e4.5}), indeed equals $0$,
we construct an involution, $\ph$ say, on this set of families that
maps a family $( T_1, T_2,\dots, T_n)$ with
associated permutation $\si$ to a family $\overline{\mathbf
T}=(\overline{ T}_1,\overline{ T}_2,\dots,\overline{ T}_n)$ with
associated permutation $\overline\si$, such that
\begin{equation}
\sgn\si=-\sgn\overline\si,\label{e4.8}
\end{equation}
and such that
\begin{equation}
 \vert{\mathbf T}\vert=\vert\overline{\mathbf T}\vert.  
\label{e4.9}
\end{equation} 
Clearly, this implies that the contribution to \eqref{e4.2} of families that
are mapped to each other cancels.

The definition of the involution $\ph$ can be copied from
\cite[Sec.~3, proof of Theorem~4]{KratBE}. For convenience, we repeat
it here.  Let $({\mathbf T},\si)$ be a pair under consideration for
the sum \eqref{e4.2}.  Besides, we assume that ${\mathbf T}$ has a crossing
point.  Consider all crossing points of two-rowed arrays with
consecutive indices (see \eqref{e4.5}).  Among these points choose those with
maximal $x$-coordinate, and among all those choose the crossing point
with maximal $y$-coordinate. Denote this crossing point by $S$. Let
$i$ be minimal such that $S$ is a crossing point of $ T_i$ and
$ T_{i+1}$. Let $ T_i$ and $ T_{i+1}$ be given by
\eqref{e4.5a} and \eqref{e4.5b}, respectively. By \eqref{e4.5}, $S$ being a crossing point
of $ T_i$ and $ T_{i+1}$ means that there exist $I$ and
$J$ such that $ T_i$ looks like
\begin{equation}
\begin{array}{rcccccl}
\tilde A^{(\si (i))}_1\le&\dots&a_{I-1}&a_I&\dots&a_{k_i}&\le \tilde
E^{(i)}_1\\
\tilde A^{(\si (i))}_2\le&\dots&b_{I-1}&b_I&\dots&b_{k_i}&\le \tilde
E^{(i)}_2,
\end{array}\label{e4.10}
\end{equation}
$ T_{i+1}$ looks like
\begin{equation}
\begin{array}{rccccccl}
\tilde A^{(\si (i+1))}_1\le&\hdotsfor2
&c_{J}&c_{J+1}&\dots&c_{k_{i+1}}&\le \tilde E^{(i)}_1\\
\tilde A^{(\si (i+1))}_2\le&\dots&d_{J-1}&d_J&\hdotsfor2
&d_{k_{i+1}}&\le \tilde E^{(i)}_2,
\end{array}\label{e4.11}
\end{equation}
$S=(a_I,d_J)$, 
{\refstepcounter{equation}\label{e4.12}}
\alphaeqn
\begin{gather}  c_J<a_I\label{e4.12a}\\
b_{I-1}<d_J\label{e4.12b}
\end{gather}
and
\begin{equation}
1\le I\le k_i+1,\quad  0\le J\le k_{i+1}.\label{e4.12c}
\end{equation}
\reseteqn
Because of the construction of $S$, the indices $I$ {\it and} $J$ are
maximal with respect to \eqref{e4.12}.

We map $({\mathbf T},\si)$ to the pair $(\overline{\mathbf
T},\si\circ(i,i+1))$ ($(i,i+1)$ denotes the transposition exchanging
$i$ and $i+1$), where $\overline{\mathbf T}=(\overline T_1,\overline
T_2,\dots,\overline T_n)$, with $\overline T_j=T_j$ for all $j\ne i,i+1$,
with $\overline T_i$ being given by
{\refstepcounter{equation}\label{e4.13}}
\alphaeqn
\begin{equation}
\begin{array} {ccccc}
\dots&c_J&a_I&\dots&a_{k_i}\phantom{,}\\
\dots&d_{J-1}&b_I&\dots&b_{k_i},
\end{array}\label{e4.13a}
\end{equation}
and with $\overline T_{i+1}$ being given by
\begin{equation}
\begin{array} {cccccc}
\hdotsfor2 &a_{I-1}&c_{J+1}&\dots&c_{k_{i+1}}\phantom{.}\\
\dots&b_{I-1}&d_J&\hdotsfor2 &d_{k_{i+1}}.
\end{array}\label{e4.13b}
\end{equation}
\reseteqn

\subsection{The properties of the involution} \label{Sub:Prop}
What we have to prove is that this operation is well-defined, i.e.,
that all the rows in \eqref{e4.13a} and \eqref{e4.13b} are strictly increasing, that
$\overline T_i$ is of type $(\si\circ(i,i+1))(i)-i= \si(i+1)-i$, that
$\overline T_{i+1}$ is of type $(\si\circ(i,i+1))(i+1)-i-1=
\si(i)-i-1$, that the bounds for the entries of $\overline T_i$ are
given by
\begin{equation*}
\begin{array} {rcccccl}
\tilde A^{(\si (i+1))}_1 \le&\dots&c_J&a_I&\dots&a_{k_i}&\le \tilde E^{(i)}_1\\
\tilde A^{(\si (i+1))}_2 \le&\dots&d_{J-1}&b_I&\dots&b_{k_i}&\le \tilde E^{(i)}_2,
\end{array}\label{e4.14}
\end{equation*}
that those for $\overline T_{i+1}$ are given by
\begin{equation*}
\begin{array} {rccccccl}
\tilde A^{(\si (i))}_1 \le&\hdotsfor2
&a_{I-1}&c_{J+1}&\dots&c_{k_{i+1}}&\le \tilde E^{(i+1)}_1\\
\tilde A^{(\si (i))}_2 \le&\dots&b_{I-1}&d_J&\hdotsfor2
&d_{k_{i+1}}&\le \tilde E^{(i+1)}_2,
\end{array}\label{e4.15}
\end{equation*}
and that \eqref{e4.4} is satisfied for $\overline T_i$ and $\overline
T_{i+1}$.
Furthermore we have to prove that $\ph$ is indeed an involution (for
which it suffices to show that \eqref{e4.12} also holds for $\overline T_i$
and $\overline T_{i+1}$), and finally we must prove \eqref{e4.8} (with
$\overline\si=\si\circ (i,i+1)$) and \eqref{e4.9}.

The claim that \eqref{e4.8} and \eqref{e4.9} hold is trivial. All other claims,
except for the claim about \eqref{e4.4}, can be proved by copying the
according arguments from the proof of Theorem~4 in \cite{KratBE} (see
the paragraphs after \cite[Eq.~(27)]{KratBE}). 

\subsection{The involution respects the ladder condition}
\label{Sub:Ladder}
It remains to show that
\eqref{e4.4} is satisfied for $\overline T_i$ and $\overline
T_{i+1}$. Unfortunately, it is necessary to supplement and refine 
the according arguments in the
proof of the main theorem in \cite{KrPrAA} (see the proof of
(4.27) and (4.28) in \cite[pp.~1035--37]{KrPrAA}) substantially in
order to cope with the situation that we encounter here.
Besides, we use the opportunity to correct an inaccuracy in
\cite{KrPrAA}.

\medskip
We have to prove that for $1\le r\le i-1$ we have 
\begin{equation} \label{e4.16a}
d_{J-i+r} < f(a_{I-1+r}),
\end{equation}
provided both $a_{I-1+r}$ and $d_{J-i+r}$ exist (if either 
$a_{I-1+r}$ or $d_{J-i+r}$ does not exist there is nothing to show), 
and 
\begin{equation}
b_{I-i+r} < f(c_{J+r}),\label{e4.16b}
\end{equation}
provided both $b_{I-i+r}$ and $c_{J+r}$ exist 
(if either $b_{I-i+r}$ or $c_{J+r}$ does not exist there is nothing
to show).

{\it Proof of \eqref{e4.16a}.}
In the following, let $r$ be fixed.  We distinguish
between two cases. If $ E_1^{(1)}\le a_{I}$, 
then we have the following chain
of inequalities:
\begin{multline} \label{e4.16c} 
d_{J-i+r}\le d_{J-1}+1-i+r
\le b_I-i+r
\le b_{I-1+r}-i+1\\
\le\tilde E_2^{(i)}-i+1
= E_2^{(i)}
\le E_2^{(1)}
<f( E_1^{(1)})
\le f(a_{I})
\le f(a_{I-1+r}),
\end{multline}
as required. (The second inequality in \eqref{e4.16c} 
follows from the fact that the
rows in \eqref{e4.13a} are strictly increasing.)

Otherwise, if $ E_1^{(1)}>a_{I}$, let us assume for the purpose
of contradiction that \eqref{e4.16a} does not hold. Then, because of
the first two inequalities in \eqref{e4.16c} we have $d_{J-i+r}\le
b_I$, and hence
\begin{equation} \label{e4.ab}
f(a_I)\le f(a_{I-1+r})\le d_{J-i+r}\le b_I.
\end{equation}
In more colloquial terms, 
the point $(a_I,b_I)$ lies outside the ladder region $L$
defined by \eqref{e3.1}. 

For the following, we make the 
conventional definitions 
$a^{(j)}_{-\si(j)+j}:=\break\tilde A^{(\si(j))}_1-1$,
$a^{(j)}_{k_j+1}:=\tilde E^{(j)}_1+1$ (which is in accordance with the
conventional definition for $a_{k+1}$ in \eqref{e4.5}), and
$b^{(j)}_0:=\tilde A^{(\si(j))}_2-1$ (which is in accordance with the
conventional definition for $b_0$ in \eqref{e4.5}).

For any $j<i$ we claim that, if for the two-rowed array $T_{j+1}$ (given by
\eqref{e4.3} with $i$ replaced by $j+1$)
we find a pair $(a^{(j+1)}_{s_{j+1}},b^{(j+1)}_{s_{j+1}})$ of entries
(i.e., $a^{(j+1)}_{s_{j+1}}$ and $b^{(j+1)}_{s_{j+1}}$ exist in
$T_{j+1}$ or are defined by means of one of the above conventional
definitions)
such that\footnote{It is at the corresponding place where the inaccuracy in
  \cite{KrPrAA} occurs. On p.~1036 the inequality chain $a_I\ge
  x_s\ge\dots\ge u_t$ has to be replaced by $a_I\ge x_s$, \dots,
  $a_i\ge u_t$, and the inequality chain $b_I\le y_s\le\dots\le v_t$
  has to be replaced by $b_I\le y_s$, \dots, $b_I\le v_t$.}
\begin{equation} \label{e4.aIbIj}
a_I\ge a^{(j+1)}_{s_{j+1}} \quad \text {and}
\quad b_I\le b^{(j+1)}_{s_{j+1}},
\end{equation}
then we can find an $h\le j$ such that the two-rowed array $T_h$ contains
a pair $(a^{(h)}_{s_h},b^{(h)}_{s_h})$ satisfying the same condition,
that is
\begin{equation} \label{e4.aIbIh}
a_I\ge a^{(h)}_{s_h} \quad \text {and}
\quad b_I\le b^{(h)}_{s_h}.
\end{equation}
In other words, we claim that if in $T_{j+1}$ we find a pair of entries
which, when considered as a lattice point, is located (weakly)
northwest of $(a_I,b_I)$, then we will also find such a pair in $T_h$
for some $h\le j$.

Let us for the moment assume that we have already established the
claim. Clearly, for $j=i-1$ the condition \eqref{e4.aIbIj} is
satisfied with $s_{j+1}=I$, in which case we have
$a^{(j+1)}_{s_{j+1}}=a^{(i)}_{I}=a_I$ and
$b^{(j+1)}_{s_{j+1}}=b^{(i)}_{I}=b_I$. Then, by iterating the
assertion of our claim, we will find that \eqref{e4.aIbIh} is
satisfied for $h=1$ and some $s_1$. 
Using this and \eqref{e4.ab} we obtain
\begin{equation*}
f(a^{(1)}_{s_1})\le f(a_I)\le b_I\le b^{(1)}_{s_1}.
\end{equation*}
However, this inequality contradicts the fact that $T_1$ obeys the
ladder condition \eqref{e4.4} 
with $i=1$ and $s=s_1$. Hence, inequality \eqref{e4.16a} must be
actually true.

\medskip
For the proof of the claim, we distinguish between four cases:

\medskip

\hbox{\hbox to 30pt{\hss(i)}
$\si(j)\ge j$ and $a^{(j)}_1\le a_I$;}

\hbox{\hbox to 30pt{\hss(ii)}
$\si(j)< j$ and $a^{(j)}_{-\si(j)+j+1}\le a_I$;}

\hbox{\hbox to 30pt{\hss(iii)}
$\si(j)\ge j$ and $a^{(j)}_1> a_I$;}

\hbox{\hbox to 30pt{\hss(iv)}
$\si(j)< j$ and $a^{(j)}_{-\si(j)+j+1}> a_I$.}

\subsubsection{Case $\si(j)\ge j$ and $a^{(j)}_1\le a_I$.}
Because we are assuming $E^{(1)}_1>a_I$, we have
$a_I\le E_1^{(1)}-1= \tilde E^{(1)}_1\le \tilde E^{(j)}_1$. Therefore
it is impossible that $a^{(j)}_1=\tilde E^{(j)}_1+1$ (by one of our
conventional assignments), and hence $a^{(j)}_1$ does indeed exist,
i.e., $k_j\ge1$ (cf\@. \eqref{e4.3} with $i$ replaced by $j$). 

Let $s_j$ be maximal such that $a^{(j)}_{s_j}\le a_I$. By the above
we have $1\le s_j\le k_j$. Therefore $b^{(j)}_{s_j}$ exists. 
If $b^{(j)}_{s_j}<b_I$, then we have
$a^{(j+1)}_{s_{j+1}}\le a_I<   a^{(j)}_{s_j+1}$ and
$b^{(j)}_{s_j} < b_I\le b^{(j+1)}_{s_{j+1}}$. But that means that
$(a^{(j)}_{s_j+1}, b^{(j+1)}_{s_{j+1}})$ is a crossing point of
$T_j$ and $T_{j+1}$ (cf\@. \eqref{e4.5c}--\eqref{e4.5e}) with
larger $x$-coordinate than 
$(a_I,d_J)$, contradicting the maximality of the crossing point
$(a_I,d_J)$. Hence, we
actually have $b^{(j)}_{s_j}\ge b_I$, and thus \eqref{e4.aIbIh} holds
with $h=j$ and with $s_j$ as above.

\subsubsection{Case $\si(j)< j$ and $a^{(j)}_{-\si(j)+j+1}\le a_I$.}
The arguments from the above case apply verbatim if one replaces
$a^{(j)}_1$ by $a^{(j)}_{-\si(j)+j+1}$ everywhere.

\subsubsection{Case $\si(j)\ge j$ and $a^{(j)}_1> a_I$.}
We show that this case actually cannot occur. Because of
\eqref{e3.5c}, we have $f(A^{(1)}_1)\le f(a_I)$, and therefore
\begin{equation*}
  b^{(j)}_0=  \tilde A^{(\si(j))}_2-1\le\tilde
  A^{(1)}_2-1=A^{(1)}_2<f(A^{(1)}_1)\le f(a_I)\le b_I\le
b^{(j+1)}_{s_{j+1}},
\end{equation*}
the two last inequalities being due to \eqref{e4.ab} and
\eqref{e4.aIbIj}.
On the other hand, we have
$a^{(j+1)}_{s_{j+1}}\le a_I<   a^{(j)}_1$. This means that 
$(a^{(j)}_1,b^{(j+1)}_{s_{j+1}})$ is a crossing point of $T_j$
and $T_{j+1}$ with larger
$x$-coordinate than $(a_I,d_J)$, which contradicts again the
maximality of $(a_I,d_J)$.

\subsubsection{Case $\si(j)< j$ and $a^{(j)}_{-\si(j)+j+1}> a_I$.}
If $b^{(j)}_{-\si(j)+j}< b_I$, then we have
$  a^{(j+1)}_{s_{j+1}}\le a_I<   a^{(j)}_{-\si(j)+j+1}$ and
$  b^{(j)}_{-\si(j)+j} < b_I\le b^{(j+1)}_{s_{j+1}}$. 
This means that 
$(a^{(j)}_{-\si(j)+j+1},b^{(j+1)}_{s_{j+1}})$ is a crossing point of $T_j$
and $T_{j+1}$ with larger
$x$-coordinate than $(a_I,d_J)$, a contradiction.
Therefore we actually have $b^{(j)}_{-\si(j)+j}\ge b_I$.

If $a^{(j)}_{-\si(j)+j}=\tilde A^{(\si(j))}_1-1\le a_I$ then \eqref{e4.aIbIh}
is satisfied with $h=j$ and $s_j=-\si(j)+j$. 
If, on the other hand, $a^{(j)}_{-\si(j)+j}>a_I$, then
of course \eqref{e4.aIbIh} cannot be satisfied for $h=j$ and any
legal $s_j$.
However, we can show that it is satisfied for some smaller $h$.

Let us pause for a moment and summarize the conditions that we are
encountering in the current case:
\begin{equation}
  \label{e4.aIbIS}
  \si(j)<j\text{, } a^{(j)}_{-\si(j)+j}>a_I \text{ and } b^{(j)}_{-\si(j)+j}\ge b_I.
\end{equation}
Clearly, there is a maximal $s$ with $s\le\si(j)\le\si(s)$. We are
going to show that we can either find an $h\le j$ and a legal $s_h$
such that \eqref{e4.aIbIh} is satisfied, or we find an index $\ell<j$
such that \eqref{e4.aIbIS} is satisfied with $j$ replaced by $\ell$
(in which case we repeat the subsequent arguments), or we can
construct a sequence of pairs $(a^{(\ell)}_{r_\ell},b^{(\ell)}_{r_\ell})$,
$r_\ell\in\{1,2,\dots ,k_\ell\}$ for $\ell\in\{s+1,s+2,\dots,j-1\}$ that satisfy
\begin{equation}\label{e4.aIbIc}
  a^{(\ell+1)}_{r_{\ell+1}}\ge a^{(\ell)}_{r_\ell}>a_I \text{ and } 
  b^{(\ell)}_{r_\ell}\ge b^{(\ell+1)}_{r_{\ell+1}}\ge b_I,
\end{equation}
where, in order that \eqref{e4.aIbIc} makes sense for $\ell=j-1$, we
set $r_j=-\si(j)+j$. 

However, if we have found such pairs for $\ell\in\{s+1,s+2,\dots,j-1\}$,
then we have
\begin{multline*} 
a_I<a^{(s+1)}_{r_{s+1}}<a^{(j)}_{-\si(j)+j}+1=\tilde A^{(\si(j))}_1\\
\le\tilde A^{(\si(j))}_1+\si(j)-s
\le\tilde A^{(\si(s))}_1+\si(s)-s\le a^{(s)}_1
\end{multline*} 
and
$$b^{(s)}_0=\tilde A^{(\si(s))}_2-1\le\tilde
A^{(\si(j))}_2-1<b^{(j)}_{-\si(j)+j}=b^{(j)}_{r_j}\le
b^{(s+1)}_{r_{s+1}}.
$$
This means that $(a^{(s)}_1,b^{(s+1)}_{r_{s+1}})$ is a crossing point
of $T_s$ and $T_{s+1}$ with larger $x$-coordinate than
$(a_I,d_J)$, contradicting again the maximality of $(a_I,d_J)$. 
Therefore we will actually find an $h\le j$ such that \eqref{e4.aIbIh}
is satisfied.

We prove our claim in \eqref{e4.aIbIc} by a reverse induction on
$\ell$. (The last two inequalities in \eqref{e4.aIbIS} guarantee that
the induction can be started.)
Suppose that we have already found indices $r_{j-1},r_{j-2},\dots,r_{\ell+1}$
satisfying \eqref{e4.aIbIc}. Then we distinguish between the two cases
$\si(\ell)\ge\ell$ and $\si(\ell)<\ell$.

First let us consider the case $\si(\ell)\ge \ell$. If
$a^{(\ell)}_1 > a^{(\ell+1)}_{r_{\ell+1}}$, then we have
$a_I < a^{(\ell+1)}_{r_{\ell+1}} < a^{(\ell)}_1$ and, if in addition $\ell\ge
\si(j)$, we have
$$b^{(\ell)}_0 = \tilde A^{(\si(\ell))}_2 -1\le \tilde
A^{(\si(j))}_2-1 <
b^{(j)}_{-\si(j)+j}=  b^{(j)}_{r_j}\le b^{(\ell+1)}_{r_{\ell+1}},
$$
where the first inequality is due to $\si(\ell)\ge\ell\ge\si(j)$.
This means that $(a^{(\ell)}_1,
b^{(\ell+1)}_{r_{\ell+1}})$ is a crossing point of $T_\ell$ and
$T_{\ell+1}$ with larger $x$-coordinate than $(a_I,d_J)$, again a
contradiction.

If $\ell<\si(j)$, we can also prove that $b^{(\ell)}_0<
b^{(\ell+1)}_{r_{\ell+1}}$, giving the same contradiction. However,
this time we must argue differently. Since $\ell> s$, all of
$\si(\ell+1),\si(\ell+2),\dots,\si(\si(j))$ must be less than $\si(j)$.
For that reason, because of $\si(\ell)\ge\ell$ 
and the pigeon hole principle, there must be a
$t\in\{\ell+1,\ell+2,\dots,\si(j)\}$ with $\si(t)<\si(\ell)$.
Then, by \eqref{e4.aIbIc}, we obtain
$$  b^{(\ell)}_0 = \tilde A^{(\si(\ell))}_2 -1\le \tilde A^{(\si(t))}_2
-1<  b^{(t)}_{r_t}\le b^{(\ell+1)}_{r_{\ell+1}}.
$$
Hence, we actually have $a^{(\ell)}_1 \le a^{(\ell+1)}_{r_{\ell+1}}$.

We also have $\tilde E^{(s)}_1
\ge\tilde A^{(\si(s))}_1+\si(s)-s-1$, because otherwise 
there would not be any two-rowed array $T_s$ (see
\eqref{e4.3} with $i=s$), i.e., the family $\mathbf T$ of two-rowed
arrays that we are considering would not exist, which is absurd.
This implies the inequality chain
\begin{multline*}
\tilde E^{(\ell)}_1
\ge\tilde E^{(s)}_1
\ge\tilde A^{(\si(s))}_1+\si(s)-s-1\\
\ge\tilde A^{(\si(j))}_1+\si(j)-s-1
\ge\tilde A^{(\si(j))}_1-1
  =a^{(j)}_{-\si(j)+j}
\ge a^{(\ell+1)}_{r_{\ell+1}}.
\end{multline*}
Therefore it is impossible that $a^{(\ell)}_1=\tilde E^{(\ell)}_1+1$
(by one of our conventional assignments), 
and hence $a^{(\ell)}_1$ does indeed exist, i.e., $k_\ell\ge1$.

Now let $r_\ell$ be maximal, such that $a^{(\ell)}_{r_\ell}\le
a^{(\ell+1)}_{r_{\ell+1}}$. By the above we have $1\le r_\ell\le
k_\ell$.
If $b^{(\ell)}_{r_\ell}< b^{(\ell+1)}_{r_{\ell+1}}$, then we have
$  a_I < a^{(\ell+1)}_{r_{\ell+1}} < a^{(\ell)}_{r_\ell+1}$ and
$  b^{(\ell)}_{r_\ell} < b^{(\ell+1)}_{r_{\ell+1}}$. This means that
$(a^{(\ell)}_{r_\ell+1},b^{(\ell+1)}_{r_{\ell+1}})$ is a crossing
point of $T_\ell$ and $T_{\ell+1}$ with larger $x$-coordinate than
$(a_I,d_J)$, which is once more a contradiction.

Hence, we actually have $b^{(\ell)}_{r_\ell}\ge
b^{(\ell+1)}_{r_{\ell+1}}$. 
Therefore, if $a^{(\ell)}_{r_\ell}\le a_I$ then \eqref{e4.aIbIh}
is satisfied with $h=\ell$ and $s_h=r_\ell$, and otherwise, 
if $a^{(\ell)}_{r_\ell}> a_I$ then \eqref{e4.aIbIc} is
satisfied.

As a last subcase, we must consider $\si(\ell)< \ell$. 
Again we have to distinguish between two cases: if
$a^{(\ell)}_{-\si(\ell)+\ell+1}\le
a^{(\ell+1)}_{r_{\ell+1}}$, we argue exactly as in the above case
where $\si(\ell)\ge \ell$ and $a^{(\ell)}_1\le a^{(\ell+1)}_{r_{\ell+1}}$. 
(We just have to replace $a^{(\ell)}_1$ by
$a^{(\ell)}_{-\si(\ell)+\ell+1}$ there.) Otherwise, 
if $a^{(\ell)}_{-\si(\ell)+\ell+1}>a^{(\ell+1)}_{r_{\ell+1}}$, we get
$b^{(\ell)}_{-\si(\ell)+\ell}\ge b^{(\ell+1)}_{r_{\ell+1}}$, because otherwise
$  a_I < a^{(\ell+1)}_{r_{\ell+1}} < a^{(\ell)}_{-\si(\ell)+\ell+1}$ and
$  b^{(\ell)}_{-\si(\ell)+\ell} < b^{(\ell+1)}_{r_{\ell+1}}$, and thus 
$(a^{(\ell)}_{-\si(\ell)+\ell+1},
b^{(\ell+1)}_{r_{\ell+1}})$ is a crossing point with
larger $x$-coordinate than $(a_I,d_J)$, again a contradiction.

Now, if $a^{(\ell)}_{-\si(\ell)+\ell}\le a_I$ then \eqref{e4.aIbIh} is satisfied
with $h=\ell$ and $s_h={-\si(\ell)+\ell}$. On the other hand, if 
$a^{(\ell)}_{-\si(\ell)+\ell}> a_I$ then 
\eqref{e4.aIbIS} is satisfied with $j$ replaced by $\ell$. In
addition we have $\ell<j$.
Consequently, we repeat the arguments subsequent to \eqref{e4.aIbIS}
with $j$ replaced by $\ell$. In that manner, we may possibly perform
several such iterations. However, these iterations must come to an
end because $\si(1)\ge 1$, and, hence, the conditions \eqref{e4.aIbIS}
cannot be satisfied for $j=1$.

\medskip
{\it Proof of \eqref{e4.16b}.}
We proceed similarly. We first
observe that we must have $a_I\le c_{J+1}$, because otherwise we
would have $c_{J+1}<a_I$ and by \eqref{e4.5d} also
$b_{I-1}<d_J<d_{J+1}$, which means that $(a_I,d_{J+1})$ is a crossing
point of $T_{i}$ and $T_{i+1}$, contradicting the
maximality of $(a_I,d_J)$. Now we distinguish again between the same two
cases as in the proof of \eqref{e4.16a}. If $ E_1^{(1)}\le a_{I}$, 
then we have the following chain of inequalities:
\begin{multline} \label{e4.16d}
b_{I-i+r}\le b_{I-1}+1-i+r
\le d_J-i+r
\le d_{J+r}-i\\
\le\tilde E_2^{(i+1)}-i
= E_2^{(i+1)}
\le E_2^{(1)}
<f( E_1^{(1)})
\le f(a_{I})
\le f(c_{J+1})
\le f(c_{J+r}),
\end{multline}
as required. (The second inequality in \eqref{e4.16d} 
follows from the fact that the
rows in \eqref{e4.13b} are strictly increasing.)
If on the other hand we have $ E_1^{(1)}>a_{I}$, then
let us assume for the purpose
of contradiction that \eqref{e4.16b} does not hold. This implies
\begin{equation*} 
f(a_I)\le f(c_{J+r})\le b_{I-i+r}< b_I.
\end{equation*}
Again, this simply means that 
the point $(a_I,b_I)$ lies outside the ladder region $L$
defined by \eqref{e3.1}. We are thus in the same situation as in the
above proof of \eqref{e4.16a}, which, in the long run, led to a
contradiction.

\bigskip
This completes the proof of the theorem.

\section{Enumeration of two-rowed arrays}
\label{Sec:5}

The entries in
the determinant in \eqref{e3.5} and \eqref{e3.6} are all generating functions $\sum
_{} ^{}z^{\vert T\vert/2}$ for two-rowed arrays $T$.  Hence, we have
to say how these can be computed. Of course, a ``nice" formula cannot
be expected in general.  There are only two cases in which ``nice"
formulas exist, the case of the trivial ladder (i.e\@., $f(x)\equiv
b+1$; see \eqref{e5.4a}), and the case of a ladder determined by a diagonal
boundary (i.e\@., $f(x)=x+D+1$, for some positive integer $D$; see
\eqref{e5.5a}). In all other cases one has to be satisfied with answers of
recursive nature.

We will describe two approaches to attack this problem. The first
leads to an extension of a formula due to Kulkarni \cite{KulkAD} (see
also \cite[Prop.~4]{KrPrAA}) for the generating function of lattice
paths with given starting and end points in a one-sided ladder
region. The second extends the alternative to Kulkarni's
formula that was proposed in \cite[Prop.~5--7]{KrPrAA}. The first approach
has the advantage of producing a formula (see Proposition~\ref{prop:4} below)
that can be compactly stated. The second approach is always at least
as efficient as the first, but is by far superior for ladder regions of a
particular kind. This is discussed in more detail 
after the proof of Proposition~\ref{prop:7}.

\begin{proposition}
\label{prop:4} 
Let $f$ be a weakly increasing function
$f:[0,a]\to [1,b+1]$ corresponding to a ladder region $L$ by means of
{\em \eqref{e3.1}}, as before. Extend $f$ to all integers by setting
$f(x):=\al_2$ for $x<0$ and $f(x):=\ep_2+1$ for $x>a$. Let
$\al_1-1<s_{k-1}<\dots<s_1<\ep_1$ be a partition of the (integer)
interval $[\al_1-1,\ep_1]$ such that $f$ is constant on each subinterval
$[s_i+1,s_{i-1}]$, $i=k,k-1,\dots,1$, with $s_k:=\al_1-1$ and
$s_0:=\ep_1$. Then the generating function $\sum _{} ^{}z^{\vert
T\vert/2}$ for all two-rowed arrays $T$ of the form {\em\eqref{e3.4}} and
satisfying {\em\eqref{e3.3}} is given by
\begin{multline} \label{e5.1}
\GF(\TA(l;(\al_1,\al_2),(\ep_1,\ep_2);f,d);z^{\vert.\vert/2})\\
=\underset {e_k-f_k=l}{\sum _{\mathbf e+d\ge\mathbf f\ge \mathbf 0} ^{}}
z^{e_k}\prod _{i=1} ^{k}\binom {s_{i-1}-s_i}{e_i-e_{i-1}} \binom
{f'(s_{i-1})-f'(s_{i})}{f_i-f_{i-1}},
\end{multline}
where $\mathbf e=(e_1,e_2,\dots,e_k)$ and $\mathbf
f=(f_1,f_2,\dots,f_k)$, where, by definition,  
$e_0=f_0=0$, where
$\mathbf e+d\ge\mathbf f\ge 0$ means $e_i+d\ge f_i\ge 0$,
$i=1,2,\dots,k$, and where $f'(x)$ agrees with $f(x)$ for $\al_1\le
x<\ep_1$, but where $f'(\al_1-1)=\al_2$ and $f'(\ep_1)=\ep_2+1$.
{\em(}All other values of $f'$ are not needed 
for the formula {\em\eqref{e5.1})}.
\end{proposition}
\begin{proof}
Let $T$ be a two-rowed array in
$\TA(l;(\al_1,\al_2),(\ep_1,\ep_2);f,d)$, represented as in \eqref{e3.4}.
Suppose that there are $e_i$ entries in the first row of $T$ that are
larger than $s_i$, and that there are $f_i$ entries in the second row
of $T$ that are larger than or equal to $f(s_i)$, $i=1,2,\dots,k$. 
Equivalently, we have
\begin{multline} \label{e5.2}
\ep_1=s_0\ge a_1>\dots>a_{e_1}>s_1\ge a_{e_1+1}>\dots>a_{e_2}>s_2\\
\ge\dots>s_{k-1}\ge a_{e_{k-1}+1}>\dots>a_{e_k}>s_k=\al_1-1,
\end{multline}
and
\begin{multline} \label{e5.3}
f'(s_0)=\ep_2+1>b_1>\dots>b_{f_1}\ge f(s_1)>b_{f_1+1}>\dots>b_{f_2}\ge f(s_2)\\
>\dots\ge f(s_{k-1})>b_{f_{k-1}+1}>\dots>b_{f_k}\ge f'(s_{k})=\al_2.
\end{multline}
In particular, we have $e_k-f_k=l$.
From \eqref{e3.3} 
it is immediate that we must have $e_i+d\ge f_i\ge 0$. Conversely,
given integer vectors $\mathbf e$ and $\mathbf f$ with $e_i+d\ge f_i\ge 0$
and $e_k-f_k=l$, by \eqref{e5.2} and \eqref{e5.3} there are
$$\prod _{i=1} ^{k}\binom {s_{i-1}-s_i}{e_i-e_{i-1}} \binom
{f(s_{i-1})-f(s_{i})}{f_i-f_{i-1}}$$ 
possible choices for the entries $a_i$ and $b_i$, $i=1,2,\dots$, in
the first and second row of a two-rowed array which satisfies
\eqref{e5.2} and \eqref{e5.3}, and thus \eqref{e3.3}.
This establishes \eqref{e5.1}. \quad \quad \qed
\end{proof}

\begin{remark}If in Proposition~\ref{prop:4} we set $l=d=0$, then 
we recover Kulkarni's formula \cite[Theorem~4]{KulkAD} (see also
\cite[Prop.~4]{KrPrAA}), because the two-rowed arrays in 
$\TA(0;(\al_1,\al_2),(\ep_1,\ep_2);f,0)$ can be interpreted as
lattice paths with starting point $(\al_1,\al_2-1)$ and end point
$(\ep_1+1,\ep_2)$ which stay in the ladder region defined by $f$.
\end{remark}

Now we describe the announced alternative method to compute the
generating function $\sum_{}^{}z^{\vert{T}\vert/2}$ for two-rowed arrays $T$
of the form \eqref{e3.4} which satisfy \eqref{e3.3}.
For sake of convenience, for $A=(\al_1,\al_2)$ and $E=(\ep_1,\ep_2)$
as before, $\al_1\le \ep_1$, we introduce the set
\begin{equation} \label{eTA*}
\TA^*(l;A,E;f,d)=\TA(l;A,E;f,d)\setminus\TA(l;A+(1,0),E;f,d),
\end{equation}
which is simply the set of those two-rowed arrays of the given form whose
first entry in the first row equals $\al_1$. 

This second method is based on the simple facts that are summarized
in Propositions~\ref{prop:5}--\ref{prop:7}. The propositions extend
in turn Propositions~5--7 in \cite{KrPrAA}.
In the following, all
binomial coefficients $\binom{n}{k}$ are understood to be equal to
zero if $n$ is negative and $k$ is positive.

\begin{proposition}
\label{prop:5}
Let $L$ be the trivial ladder determined by
the function $f(x)\equiv b+1$ by means of {\em \eqref{e3.1}}. Let $A=(\al_1,\al_2)$
and $E=(\ep_1,\ep_2)$ be lattice points and $l$ and $d$ arbitrary
integers. Then we have
\begin{equation}
\GF\big(\TA(l;A,E;f,d);z^{\vert.\vert/2}\big)
=\sum_k\binom{\ep_1-\al_1+1}{k+l}
       \binom{\ep_2-\al_2+1}{k}z^{k+{l}/{2}},
\label{e5.4a}
\end{equation}
and if $\al_1\le \ep_1$ we have 
\begin{equation}
\GF\big(\TA^*(l;A,E;f,d);z^{\vert.\vert/2}\big)
=\sum_k\binom{\ep_1-\al_1}{k+l-1}
       \binom{\ep_2-\al_2+1}{k}z^{k+{l}/{2}}.
\label{e5.4b}
\end{equation} 
\end{proposition}

\begin{proposition}
\label{prop:6}
Let $L_D$ be a ``diagonal" ladder determined
by the function $f(x)=x+D+1$ for an integer $D$ by means of {\em
\eqref{e3.1}}. Let $d$ be a nonnegative integer and $l$ an integer such that
$l+d\ge 0$. Let $A=(\al_1,\al_2)$ and $E=(\ep_1,\ep_2)$ be lattice
points such that $\al_1+D+1+l+d\ge\al_2$ and
$\ep_1+D+1+d\ge\ep_2$. Then we have
\begin{multline} 
\GF\big(\TA(l;A,E;f,d);z^{\vert.\vert/2}\big)
=\sum_k\Bigg(\binom{\ep_1-\al_1+1}{k+l}
              \binom{\ep_2-\al_2+1}{k}\\
            -\binom{\ep_1-\al_2+D+1}{k-d-1}
              \binom{\ep_2-\al_1-D+1}{k+l+d+1}\Bigg)z^{k+{l}/{2}},
\label{e5.5a}
\end{multline}
and if $\al_1\le \ep_1$ we have
\begin{multline} 
\GF\big(\TA^*(l;A,E;f,d);z^{\vert.\vert/2}\big)
=\sum_k\Bigg(\binom{\ep_1-\al_1}{k+l-1}
              \binom{\ep_2-\al_2+1}{k}\\
            -\binom{\ep_1-\al_2+D+1}{k-d-1}
              \binom{\ep_2-\al_1-D}{k+l+d}\Bigg)z^{k+{l}/{2}}.
\label{e5.5b}
\end{multline}
\end{proposition}

\medskip
\noindent
{\it Proof of Propositions~\ref{prop:5} and \ref{prop:6}}.
Identities \eqref{e5.4a} and \eqref{e5.4b} are immediate from the
definitions. 

To prove identity \eqref{e5.6}, we note that the number of two-rowed arrays
{\refstepcounter{equation}\label{e5.6}}
\alphaeqn
\begin{equation}
\begin{array}{rcccccccl}
\al_1\le&a_{-l+1}&a_{-l+2}&\dots&a_0&a_1&\dots&a_k&\le\ep_1\\
\al_2\le&        &        &     &   &b_1&\dots&b_k&\le\ep_2
\end{array}\label{e5.6a}
\end{equation}
that obey 
\begin{equation}
b_i<a_{i+d}+D+1, \quad i=1,2,\dots,k,\label{e5.6b}
\end{equation} 
\reseteqn
is the number of
{\it all} two-rowed arrays of the form \eqref{e5.6a} minus those that violate
the condition \eqref{e5.6b}. Clearly, the generating function for the former
two-rowed arrays is given by the first term in the sum on the right
hand side of \eqref{e5.5a}. We claim that the two-rowed arrays of the form
\eqref{e5.6a} that violate \eqref{e5.6b} are in one-to-one correspondence with
two-rowed arrays of the form
\begin{equation}
\begin{array}{rcccccccl}
\al_2-D\le&           &         &     &   &c_1&\dots&c_{k-d-1}&\le\ep_1\\
\al_1+D\le&d_{-l-2d-1}&d_{-l-2d}&\dots&d_0&d_1&\dots&d_{k-d-1}&\le\ep_2.
\end{array}\label{e5.6c}
\end{equation}
(In particular, if $k\le d$ then there is no two-rowed array of the
form \eqref{e5.6c}, in agreement with the fact that there cannot be
any two-rowed array of the form \eqref{e5.6a} violating \eqref{e5.6b} 
in that case.)
The generating function for the two-rowed arrays in \eqref{e5.6c} is
$$\sum_k\binom{\ep_1-\al_2+D+1}{k-d-1}
\binom{\ep_2-\al_1-D+1}{k+l+d+1}z^{k+l/2},$$ 
which is exactly the negative of the second term on the right-hand
side of \eqref{e5.5a}. This would prove \eqref{e5.5a}.
So it remains to
construct the one-to-one correspondence.

The correspondence that we are going to describe is
gleaned from \cite{KrMoAC}, see also \cite[Sec.~13.4]{KratBL} and
\cite{KrNiAA}.
Take a two-rowed array of the form \eqref{e5.6a} that violates condition
\eqref{e5.6b}, i.e., there is an index $i$ such that $b_i\ge
a_{i+d}+D+1$. Let $I$ be the largest integer with this property. Then
map this two-rowed array to
$$\begin{array}{rccccccccl}
\al_2\!-\!D\le&&(b_1\!-\!D)&\hdotsfor2&(b_{I-1}\!-\!D)&a_{I+d+1}&\dots&a_k&\le\ep_1\\
\al_1\!+\!D\le&(a_{-l+1}\!+\!D)&\hdotsfor1&(a_{I+d}\!+\!D)&b_I&\hdotsfor3&b_k&\le\ep_2.
\end{array}$$
Note that both rows are strictly increasing because of $b_{I-1}-D\le
b_{I+1}-D-2<a_{I+d+1}$. If $I=1$, we have to check in addition that
$\al_2-D\le a_{d+2}$, which is indeed the case, because 
$$ 
a_{d+2}\ge a_{d+1}+1\ge\dots
\ge a_{-l+1}+1+l+d\ge\al_1+1+l+d\ge\al_2-D.
$$ 
Similarly, it can be checked that $b_{I-1}-D\le\ep_1$ if $I=k-d$.
It is easy to see that the array is of
the form \eqref{e5.6c}.

The inverse of this map is defined in the same way. Take a two-rowed
array of the form \eqref{e5.6c}. Let $J$ be the
largest integer such that $d_J\ge c_{J+d}+D+1$, if existent. 
If there is no such
integer, then let $J=-d$. We map this two-rowed array to
$$\begin{array}{rrcccccccl}
\al_1\le&(d_{-l-2d-1}\!-\!D)\quad&\hdotsfor3&(d_{J\kern-.5pt-\kern-.5pt1}
\!-\!D)&c_{J+d+1}
&\dots&c_{k\kern-.5pt-\kern-.5ptd\kern-.5pt-\kern-.5pt1}&\le\ep_1\\
\al_2\le&          (c_1\!+\!D)&\dots&(c_{J+d}\!+\!D)&d_J&\hdotsfor3
      &d_{k\kern-.5pt-\kern-.5ptd\kern-.5pt-\kern-.5pt1}&\le\ep_2
\end{array}$$
Since we required $l+d\ge 0$ the entry $d_{J-1}-D$ exists even if
$J=-d$. This implies that the two-rowed array we obtained violates
condition \eqref{e5.6b}, since $d_J\ge d_{J-1}+1=(d_{J-1}-D)+D+1$. As
above, it can be checked that both rows are strictly increasing, even
in the case $J=-d$, and that the array is of the correct form. 

Equation~\eqref{e5.5b} is an immediate consequence of \eqref{e5.5a}
and the definition \eqref{eTA*} of $\TA^*(l;A,E;f,d)$. 
\quad \quad \qed
\medskip

\begin{proposition}
\label{prop:7} 
Let $L$ be an arbitrary ladder  
given by a function $f$ by means of {\em\eqref{e3.1}}, 
let $A=(\al_1,\al_2)$,
$E=(\ep_1,\ep_2)$ be lattice points in $L$, and let $d$ be a
nonnegative integer and $l$ an integer such that $l+d\ge 0$. Then for
all $x\in[0,a]$ such that $\al_2\le f(x)\le \ep_2+1$ we have
\begin{align} \notag
\GF\big(\TA(l;A,E;f&,d);z^{\vert.\vert/2}\big)\\
\notag
=\sum_{j=x+1}^{\ep_1}
&\GF\big(\TA(l+d;A,(j-1,f(x)-1);f,0);
         z^{\vert.\vert/2}\big)\\
\notag
\cdot
&\GF\big(\TA^*(-d;(j,f(x)),E;f,d);
         z^{\vert.\vert/2}\big)\\
\notag
+\sum_{e=0}^d
&\GF\big(\TA(l+d-e;A,(\ep_1,f(x)-1);f,e);
         z^{\vert.\vert/2}\big)\\
\cdot
&\binom{\ep_2-f(x)+1}{d-e} z^{(d-e)/{2}}.
\label{e5.7}
\end{align}
\end{proposition}
\begin{proof} 
We show this recurrence relation by decomposing an array
\begin{equation}
\begin{array}{rccccccccl}
\al_1\le&a_{-l+1}&a_{-l+2}&\dots&a_{-1}&a_0&a_1&\dots&a_k&\le\ep_1\\
\al_2\le&        &        &     &      &   &b_1&\dots&b_k&\le\ep_2
\end{array}\label{e5.8}
\end{equation}
in $\TA(l;A,E;f,d)$ --- the generating function of which is the left-hand 
side of \eqref{e5.7} --- into two parts. Let $I$ be the smallest integer
with $b_I\ge f(x)$, or, if all $b_I$ are smaller than $f(x)$, let
$I=k+1$. Now we have to distinguish between two cases.

If $I+d<k+1$, we decompose such an array into the array
\begin{equation*}
\begin{array}{rccccccccccl}
\al_1\le&a_{-l+1}&a_{-l+2}&\dots&a_{-1}&a_0&a_1&\dots&a_{d+1}&\dots&a_{I-1+d}&\le a_{I+d}-1\\
\al_2\le&        &        &     &      &   &   &     &b_1    &\dots&b_{I-1}  &\le f(x)-1
\end{array}\label{e5.9a}
\end{equation*} 
in $\TA(l+d;A,\big(a_{I+d}-1,f(x)-1\big);f,0)$, and the array
\begin{equation*}
\begin{array}{rcccccl}
a_{I+d}\le&   &       &a_{I+d}&\dots&a_k&\le\ep_1\\
f(x)\le   &b_I&\dots  &b_{I+d}&\dots&b_k&\le\ep_2
\end{array}\label{e5.9b}
\end{equation*}
in $\TA^*(-d;\big(a_{I+d},f(x)\big),E;f,d)$. Clearly, this is a
pair of two-rowed arrays enumerated by the first sum in the right hand
side of \eqref{e5.7}, with the summation index $j$ equal to $a_{I+d}$.

If $I+d\ge k+1$, we decompose \eqref{e5.8} into the array
\begin{equation*}
\begin{array}{rccccccccccl}
\al_1\le&a_{-l+1}&a_{-l+2}&\dots&a_{-1}&a_0&a_1&\dots&a_{k-I+2}&\dots&a_k      &\le \ep_1\\
\al_2\le&        &        &     &      &   &   &     &b_1    &\dots&b_{I-1}  &\le f(x)-1
\end{array}\label{e5.10a}
\end{equation*}
in $\TA(l-I+k+1;A,\big(\ep_1,f(x)-1\big);f,d+I-k-1)$, and a single
row
\begin{equation*}
\begin{array}{rcccl}
f(x)\le &b_I&\dots  &b_k&\le\ep_2.
\end{array}\label{e5.10b}
\end{equation*}
Note that, if $I=k+1$, this row is empty. These pairs are enumerated
by the second sum on the right hand side of \eqref{e5.7}, with the summation
index $e$ equal to $d+I-k-1$.
\quad \quad \qed
\end{proof}

Now, here is the second method for determining 
$\GF(\TA(l;(\al_1,\al_2),\break(\ep_1,\ep_2);f,d));z^{\vert.\vert/2})$ 
for any given ladder $L$ of the form \eqref{e3.1}, with
points $A=(\al_1,\al_2)$ and $E=(\ep_1,\ep_2)$ located inside $L$:
partition the {\it border\/} of $L$, i.e\@., the set of points
$\{(x,f(x)):x\in[0,a]\}$ into horizontal and diagonal pieces, say
$L^1,L^2,\dots,L^m$, where $L^i=\{(x,f(x)):x_{i-1}<x\le x_i\}$, for
some $-1=x_0<x_1<x_2<\dots<x_m=a$, each $L^i$ being either horizontal
or diagonal. Then apply the recurrence \eqref{e5.7} in succession with
$x=x_{m-1},x_{m-2},\dots,x_1$ and use \eqref{e5.4a}--\eqref{e5.5b} to compute all
the occurring generating functions. 


To give an example, in the case of
the ladder of Figure~2 we would choose $m=3$, $x_1=3$, $x_2=7$,
$x_3=13$, and the resulting formula reads
\begin{align} \notag
\GF\big(\TA&(l;A,E;f,d);z^{\vert.\vert/2}\big)
=\sum_{j=8}^{\ep_1}\sum_{k\ge0}z^{k-{d}/{2}}
\binom{\ep_1-j}{k-d-1}\binom{\ep_2-12}{k}\\
\notag
&\kern-9pt\cdot\Bigg(\sum_{i=4}^{j-1}\sum_{k_1,k_2\ge0}
z^{k_1+k_2+(l+d)/{2}}\binom{i-\al_1}{k_1+l+d}\binom{7-\al_2}{k_1}\\
\notag
&\hphantom{\Bigg(\sum_{i=4}^{13}\sum_{k_1,k_2\ge0}z^{k_1+k_2+(l+d)/{2}}\Bigg.}
\kern-9pt\cdot
\lK\binom{j-i-1}{k_2-1}\binom{6}{k_2}-\binom{7-i}{k_2}\binom{j-2}{k_2-1}\rK\\
\notag
&+\sum_{k_1\ge0}z^{k_1+(l+d)/{2}}\binom{j-\al_1}{k_1+l+d}
\binom{7-\al_2}{k_1}\Bigg)\\
\notag
+\sum_{e=0}^d&z^{(d-e)/{2}}\binom{\ep_2-12}{d-e}\\
\notag
&\cdot\Bigg(\sum_{i=4}^{\ep_1}\sum_{k_1,k_2\ge0}z^{k_1+k_2+(l+d-e)/{2}}
\binom{i-\al_1}{k_1+l+d}\binom{7-\al_2}{k_1}\\
\notag
&\hskip3cm
\cdot
\Bigg(\binom{\ep_1-i}{k_2-e-1}\binom{6}{k_2}
-\binom{7-i}{k_2}\binom{\ep_1-1}{k_2-e-1}\Bigg)\\
&+\sum_{f=0}^{e}\sum_{k\ge0}z^{k+(l+d+e)/{2}-f}
\binom{\ep_1-\al_1+1}{k+l+d-f}\binom{7-\al_2}{k}
\binom{6}{e-f}\Bigg).
\label{e5.11}
\end{align}

If $L$ consists of not too many pieces, both methods are feasible
methods, see our Example in Section~\ref{Sec:3}. Both methods yield
$(2m-1)$-fold sums if the partition of the border consists of
horizontal pieces throughout.  However, the second method is
by far superior in case of long diagonal portions in the border of $L$,
since then Kulkarni's formula involves a lot more summations. For
example, when we implemented formula \eqref{e5.11} (in {\sl
Mathematica}) it was by a
factor of 40.000 (!) faster than the corresponding implementation of
formula \eqref{e5.1}. (Indeed, the ``simplicity" of the formula
\eqref{e5.1} in comparison to \eqref{e5.11} is deceptive, as
\eqref{e5.1} involves an 11-fold summation in that case, whereas
\eqref{e5.11} has only 3-fold, 4-fold, and 5-fold sums.) Of
course, in the worst case, when $L$ consists of $1$-point pieces
throughout, both methods are nothing else than plain counting, and
therefore useless.  For computation in case of such ``fractal"
boundaries it is more promising to avoid Theorem~\ref{thm:2} and instead try to
extend the dummy path method in \cite{KrMoAA} such that it also
applies to the enumeration of nonintersecting lattice paths with
respect to turns.

\clearpage
\addcontentsline{toc}{section}{Index}
\flushbottom
\printindex


\begin{thebibliography}{77}
%
\addcontentsline{toc}{section}{References}


\bibitem{AbhyAB} S. S. Abhyankar, 
{\it Enumerative combinatorics of Young tableaux},
Marcel Dekker, New York, Basel, 1988.

\bibitem{AbKuAC} S. S. Abhyankar and D. M. Kulkarni, 
{\it On Hilbertian ideals},
 Linear Alg\@. Appl\@. {\bf 116}
(1989), 53--76.

\bibitem{BiLaAA} 
S. C. Billey and V. Lakshmibai, {\it Singular loci of Schubert
varieties}, Birkh\"auser, Boston, 2000.

\bibitem{BrHeAA} W.    Bruns and J. Herzog, 
{\it On the computation of $a$-invariants},
 Manuscripta Math\@. {\bf 77}
(1992), 201--213.

\bibitem{ConcAB} A.    Conca, 
{\it Ladder determinantal rings}, 
 J.~Pure Appl\@. Algebra  {\bf 98} 
(1995), 119--134.

\bibitem{CoHeAA} A.    Conca and J. Herzog, 
{\it On the Hilbert function of determinantal rings and their canonical module}, 
 Proc\@. Amer\@. Math\@. Soc\@.  {\bf 122} 
(1994), 677--681.

\bibitem{GhorAC} S. R. Ghorpade, 
{\it Abhyankar's work on Young tableaux and some recent
developments}, in:
Proc\@. Conf\@. on Algebraic Geometry and Its Applications 
(Purdue Univ\@., June 1990), Sprin\-ger--Ver\-lag, New York,
1994, pp.~215--249.

\bibitem{GhorAD} S. R. Ghorpade, {\it Young bitableaux, lattice 
paths and Hilbert functions}, J. Statist\@. Plann\@. Inference {\bf 
54} (1996), pp.~55--66.

\bibitem{GhorAF} S. R. Ghorpade, 
{\it Hilbert functions of ladder determinantal varieties},
Discrete Math\@. (to appear).

\bibitem{GoLaAA} N. Gonciulea and V. Lakshmibai, {\it Singular loci
of Schubert varieties and ladder determinantal varieties}, J. Algebra
{\bf 229} (2000), 463--497.

\bibitem{HeTrAA} J.    Herzog and N. V. Trung, 
{\it Gr\"obner bases and multiplicity of determinantal and Pfaffian ideals},
 Adv\@. in Math\@. {\bf 96}
(1992), 1--37.

\bibitem{KratAF} C.    Krattenthaler, 
{\it Counting lattice paths with a linear boundary, Part 2: $q$-ballot and $q$-Catalan numbers},
 Sitz\@.ber\@. d\@. \"OAW, Math-na\-tur\-wiss\@. Klasse {\bf 198}
(1989), 171--199.

\bibitem{KratAP} C.    Krattenthaler, 
{\it The major counting of nonintersecting lattice paths and generating functions for tableaux},
 Mem\@. Amer\@. Math\@. Soc\@. 115, no.~552,
 Providence, R.~I., 1995.

\bibitem{KratBE} C.    Krattenthaler, 
{\it Counting nonintersecting lattice paths with turns},
 S\'eminaire Lotharingien Combin\@. {\bf 34}
(1995), paper B34i, 17~pp.

\bibitem{KratBL} C.    Krattenthaler, {\em The enumeration of lattice 
paths with respect to their number of turns}, in: Advances in 
Combinatorial Methods and Applications to Probability and Statistics,
N.~Balakrishnan, ed., Birkh\"auser, Boston, 1997, pp.~29--58.

\bibitem{KrNiAA} C.    Krattenthaler and H. Niederhausen, {\it Lattice 
paths with weighted left turns above a parallel to the diagonal}, 
Congr\@. Numer\@. {\bf 124} (1997), 73--80.

\bibitem{KrPrAA} C.    Krattenthaler and M. Prohaska, 
{\it A remarkable formula for counting nonintersecting lattice paths in a ladder with respect to turns},
 Trans\@. Amer\@. Math\@. Soc\@.  {\bf 351}
(1999), 1015--1042.

\bibitem{KrMoAC} C.    Krattenthaler and S. G. Mohanty, 
{\it On lattice path counting by major and descents},
 Europ\@. J. Combin\@. {\bf 14}
(1993), 43--51.

\bibitem{KrMoAA} C.    Krattenthaler and S. G. Mohanty, 
{\it Counting tableaux with row and column bounds},
 Discrete Math\@. {\bf 139}
(1995), 273--286.

\bibitem{KulkAD} D. M. Kulkarni, 
{\it Hilbert polynomial of a certain ladder-determinantal ideal},
 J. Alg\@. Combin\@. {\bf 2}
(1993), 57--72.

\bibitem{KulkAC} D. M. Kulkarni, 
{\it Counting of paths and coefficients of Hilbert polynomial of a determinantal ideal},
 Discrete Math\@. {\bf 154}
(1996), 141--151.

\bibitem{ModaAA} M. R. Modak, 
{\it Combinatorial meaning of the coefficients of a Hilbert polynomial},
 Proc\@. Indian Acad\@. Sci\@. (Math\@. Sci\@.) {\bf 102}
(1992), 93--123.

\bibitem{RubeAB} M. Rubey, 
{\it Comment on `Counting nonintersecting lattice paths with turns'
by C.~Krattenthaler},  S\'eminaire Lotharingien Combin\@. {\bf 34},
Comment on paper B34i, 2001.


\end{thebibliography}
\end{document}